\documentclass[preprint,12pt]{elsarticle}

\usepackage{amsmath} % assumes amsmath package installed
\usepackage{amssymb}  % assumes amsmath package installed
\usepackage{amsxtra}

\newtheorem{lemma}{Lemma}[section]
\newtheorem{theorem}[lemma]{Theorem}

\newtheorem{definition}[lemma]{Definition}

\newtheorem{remark}[lemma]{Remark}
\DeclareMathOperator{\sign}{sign}
\newcommand{\Simp}{\ensuremath{\mathcal{S}}}
\newcommand{\R}{\ensuremath{\mathbb{R}}}
\newcommand{\N}{\ensuremath{\mathbb{N}}}

\newcommand {\cK}      {\ensuremath{\mathcal{K}}}
\newcommand \cL      {\mathcal{L}}

\newcommand \D   {\text{D}}
\newcommand \Id  {\operatorname{Id}}
\newcommand{\K}{\ensuremath{\mathcal{K}}}
\newcommand{\KL}{\ensuremath{\mathcal{KL}}}

\newcommand \eps   {\varepsilon}

\journal{Nonlinear Analysis: Hybrid Systems}

\begin{document}

\begin{frontmatter}

%% Title, authors and addresses

%% use the tnoteref command within \title for footnotes;
%% use the tnotetext command for the associated footnote;
%% use the fnref command within \author or \address for footnotes;
%% use the fntext command for the associated footnote;
%% use the corref command within \author for corresponding author footnotes;
%% use the cortext command for the associated footnote;
%% use the ead command for the email address,
%% and the form \ead[url] for the home page:
%%
%% \title{Title\tnoteref{label1}}
%% \tnotetext[label1]{}
%% \author{Name\corref{cor1}\fnref{label2}}
%% \ead{email address}
%% \ead[url]{home page}
%% \fntext[label2]{}
%% \cortext[cor1]{}
%% \address{Address\fnref{label3}}
%% \fntext[label3]{}

\title{Stability of interconnected impulsive systems with and without time-delays using Lyapunov methods\tnoteref{t2}} 
%\tnotetext[t1]{A part of the results of this paper was presented at the 8th IFAC Symposium on Nonlinear Control Systems (NOLCOS), September, 01-03, 2010, Bologna, Italy.}
\tnotetext[t2]{Sergey Dashkovskiy, Andrii Mironchenko and Lars Naujok are supported by the German Research Foundation (DFG) as part of the Collaborative Research Center 637 ''Autonomous Cooperating Logistic Processes: A Paradigm Shift and its Limitations''. Michael Kosmykov is supported by the Volkswagen Foundation as part of the Project Nr. I/82684 ''Dynamic Large-Scale Logistics Networks''. The authors thank Professor Andrew Teel for the valuable conversations and the anonymous reviewers for the helpful suggestions.}

%% use optional labels to link authors explicitly to addresses:
%% \author[label1,label2]{<author name>}
%% \address[label1]{<address>}
%% \address[label2]{<address>}

\author[Dash]{Sergey Dashkovskiy}\ead{sergey.dashkovskiy@fh-erfurt.de}    % Add the 
\author[Brem]{Michael Kosmykov}\ead{kosmykov@math.uni-bremen.de}               % e-mail address 
\author[Brem]{Andrii Mironchenko}\ead{andmir@math.uni-bremen.de}
\author[Brem]{Lars Naujok\corref{cor1}}\ead{larsnaujok@math.uni-bremen.de}  % (ead) as shown

%\thanks[Corr]{Corresponding author. Tel.: +49 421 21863744; Fax: +49 421 2189863744}
\cortext[cor1]{Corresponding author. Tel.: +49 421 21863744; Fax: +49 421 2189863744}
\address[Dash]{University of Applied Sciences Erfurt, Department of Civil Engineering, Erfurt, Germany}
\address[Brem]{University of Bremen, Centre of Industrial Mathematics, P.O.Box 330440, 28334 Bremen, Germany}  % Please supply                                              

\begin{abstract}
In this paper, we consider input-to-state stability (ISS) of impulsive control systems with and without time-delays. We prove that if the time-delay system possesses an exponential Lyapunov-Razumikhin function or an exponential Lyapunov-Krasovskii functional, then the system is uniformly ISS provided that the average dwell-time condition is satisfied. Then, we consider large-scale networks of impulsive systems with and without time-delays and prove that the whole network is uniformly ISS under the small-gain and the average dwell-time condition. Moreover, these theorems provide us with tools to construct a Lyapunov function (for time-delay systems - a Lyapunov-Krasovskii functional or a Lyapunov-Razumikhin function) and the corresponding gains of the whole system, using the Lyapunov functions of the subsystems and the internal gains, which are linear and satisfy the small-gain condition. We illustrate the application of the main results on examples.
\end{abstract}


\end{frontmatter}

\section{Introduction}
Impulsive systems combine continuous and discontinuous behavior of a dynamical system \cite{HCN06}. The continuous dynamics is typically described by differential equations and the discontinuous behavior are instantaneous state jumps that occur at given time instants, also referred to as impulses. Impulsive systems are closely related to hybrid systems \cite{HCN06} and switched systems \cite{SWM07} and have wide range of applications. %in network control, engineering, biological or economical systems, see \cite{ScS00,HCN06, SWM07}.

In this paper we study the input-to-state stability (ISS) property of impulsive systems. ISS was first introduced for continuous systems in \cite{Son89}. Lyapunov functions provide a useful tool to verify the ISS property (see \cite{SoW95}) as well as for other variants of ISS, namely input-to-state dynamical stability \cite{Gru02a}, local ISS \cite{SoW96,DaR10} and integral ISS (iISS) \cite{Son98}. Investigation of ISS for hybrid systems can be found in \cite{CaT09}. For time-delay systems the ISS property can be verified by Lyapunov-Razumikhin functions \cite{Tee98} or Lyapunov-Krasovskii functionals \cite{PeJ06}.

For impulsive systems the ISS and iISS properties were studied in \cite{HLT08} for the delay-free case and in \cite{ChZ09} for non-autonomous time-delay systems. Sufficient conditions, which assure ISS and iISS of an impulsive system, were derived using locally Lipschitz continuous exponential ISS-Lyapunov(-Razumikhin) functions. In \cite{HLT08} the average dwell-time condition was used, whereas in \cite{ChZ09} a fixed dwell-time condition was utilized. The average dwell-time condition was introduced in \cite{HeM99} for switched systems. 

In this paper we provide a Lyapunov-Krasovskii type and a Lyapunov-Razumikhin type ISS theorem for single impulsive time-delay systems using the average dwell-time condition. The proofs use the idea of the proof of \cite{HLT08}. For the Razumikhin type ISS theorem we require as an additional condition that the Lyapunov gain fulfills a small-gain condition. To prove this theorem we show the equivalence of ISS and the conjunction of two properties, one of them was also used in ISS definition from \cite{Tee98}, namely uniform global stability and a uniform convergence property. In contrast to the Razumikhin-type theorem from \cite{ChZ09} we consider autonomous time-delay systems and the average dwell-time condition. This condition considers the average of impulses over an interval, whereas the fixed dwell-time condition considers the (minimal or maximal) interval between two impulses. Our theorem allows to verify the ISS property for larger classes of impulse time sequences, however, we have used an additional technical condition on the Lyapunov-gain in our proofs.

Considering large-scale networks of impulsive systems our main goal is to find sufficient conditions which assure ISS of interconnections of impulsive systems with and without time-delays. To this end, we use the approach used for networks of continuous systems. The first results about the ISS property for the delay-free case were given for two coupled continuous systems in \cite{JTP94} and for an arbitrarily large number ($n\in\N$) of coupled continuous systems in \cite{DRW07}, using a small-gain argument which is a condition on the interconnection structure of the system. Lyapunov versions of the ISS small-gain theorems were proved in \cite{JMW96} (two systems) and \cite{DRW10} ($n$ systems), where ISS-Lyapunov functions for the overall system are constructed. There are also known results for the ISS property of hybrid systems, see \cite{NeT08} (two systems) and \cite{DaK09} ($n$ systems), as well as for the infinite-dimensional systems \cite{DaM11}. 

In \cite{DaN10} Lyapunov-Razumikhin functions and Lyapunov-Krasovskii functionals are used to verify the ISS property of large-scale time-delay systems, where a small-gain condition is used. An approach with vector Lyapunov functions can be found in \cite{KaJ11}. 

We prove that under a small-gain condition with linear gains and a dwell-time condition according to \cite{HLT08} a large-scale network of impulsive systems has the ISS property and construct the exponential ISS-Lyapunov, Lyapunov-Razumikhin and Lyapunov-Krasovskii function(al) and the corresponding gains of the whole system. %Using exponential Lyapunov functions or functionals for interconnected systems we use linear gains to derive the dwell-time condition for the interconnection. 

The paper is organized as follows: In Section~2 we note some basic definitions. Single impulsive systems are studied in Section~3, where the Lyapunov-Krasovskii and the Lyapunov-Razumikhin methodologies for impulsive time-delay systems are introduced, including the first main results of this paper. In Section~4, large-scale networks of impulsive systems with and without time-delays are considered and three different Lyapunov-type theorems are proved. An illustrative example of the application of the main results for networks can be found in Section~5. Finally, Section~6 concludes this paper with a short summary.

%%%%%%%%%%%%%%%%%%%%%%%%%%%%%%%%%%%%%%%%%%%%%%%%%%%%%%%%%%%%%%%%%%%%%%%%%%%%%%%%%%%%%%%%%%%%%%%%%%%%%%%%%%%%%%%%%%%%%%%%%%%%%%%%%%%%%%%%%%%%%%%%%%%%%%%%%%%%%%%%%%%  
\section{Preliminaries}
By $x^T$ we denote the transposition of a vector $x\in\R^N,$ $N\in\N$, furthermore $\R_+:=[0,\infty)$ and $\R^N_+$ denotes the positive orthant $\left\{x\in\R^N:x\geq0 \right\}$ where we use the partial order for $x,y\in\R^N$ given by
\begin{align*}
&x\geq y\Leftrightarrow x_i \geq y_i,\ i=1,\ldots,N\text{ and } x\not\geq y\Leftrightarrow\exists i:x_i<y_i,\\
&x > y\Leftrightarrow x_i > y_i,\ i=1,\ldots,N.%\text{ and } x\not> y\Leftrightarrow\exists i: x_i\leq y_i.
\end{align*}
We denote the Euclidean norm by $\left|\cdot\right|$. For a piecewise continuous function $x:I \to  \R^N$ we define $\left\|x\right\|_{\text{I}}:=\sup\limits_{t\in I}|x(t)|$. $\nabla V$ denotes the gradient of a function $V$. The upper right-hand side derivative of a locally Lipschitz continuous function $V:\R^N\rightarrow\R_+$ along trajectory $x(\cdot)$ 
 is defined by 
\begin{align*}
	\D^+V(x(t))=\limsup\limits_{h \rightarrow 0^+}\frac{V(x(t+h))-V(x(t))}{h}.
\end{align*}
The function $x^t:\left[-\theta,0\right]\rightarrow\R^N$ is given by $x^t(\tau):=x(t+\tau),\ \tau\in\left[-\theta,0\right]$, where $\theta\in\R_+$ is the maximum involved delay and we denote the norm $\|x^t\|:=\max_{t-\theta\leq s \leq t}\left|x(s)\right|$. For $a,b\in\R,$ $a<b$, let $PC\left(\left[a,b\right];\R^N\right)$ denote the Banach space of piecewise right-continuous functions defined on $\left[a,b\right]$ equipped with the norm $\left\|\cdot\right\|_{\left[a,b\right]}$ and take values in $\R^N$.

To define the stability notion we use the following classes.
\begin{definition}
Classes of comparison functions are:
\begin{equation}\notag
\begin{array}{ll}
{\cK} &:= \left\{\gamma:\R_+\rightarrow\R_+\left|\ \gamma\mbox{ is continuous, }\gamma(0)=0\right.\text{and strictly increasing} \right\},\\
%{\cK} &:= \left\{\gamma:\R_+\rightarrow\R_+\left|\ \gamma\mbox{ is continuous, }\gamma(0)=0\right.\right.\\
%&\phantom{:= \left\{\gamma:\R_+\rightarrow\R_+\ \ \right. }\left.\text{and strictly increasing}\right\},\\
{\cK_{\infty}}&:=\left\{\gamma\in\cK\left|\ \gamma\mbox{ is unbounded}\right.\right\},\\
{\cL}&:=\left\{\gamma:\R_+\rightarrow\R_+\left|\ \gamma\mbox{ is continuous and decreasing}\right.\right.\\
&\phantom{:=\left\{\gamma:\R_+\rightarrow\R_+\left|\ \right.\right.} \text{ with } \lim_{t\rightarrow\infty}\gamma(t)=0\},\\
{\cK\cL} &:= \left\{\beta:\R_+\times\R_+\rightarrow\R_+\left|\ \right. \beta(\cdot,t)\in{\cK},\ \beta(r,\cdot)\in {\cL},\ \forall t,r \geq 0\right\}.
%{\cK\cL} &:= \left\{\beta:\R_+\times\R_+\rightarrow\R_+\left|\ \beta \mbox{ is continuous,}\right.\right.\\
%&\phantom{:= \left\{\beta:\R_+\times\right.}\left.\beta(\cdot,t)\in{\cK},\ \beta(r,\cdot)\in {\cL},\ \forall t,r \geq 0\right\}.
\end{array}
\end{equation}
\end{definition}
Note that for $\gamma\in{\cK}_{\infty}$ the inverse function $\gamma^{-1}$ always exists and $\gamma^{-1}\in{\cK}_{\infty}$.
%%%%%%%%%%%%%%%%%%%%%%%%%%%%%%%%%%%%%%%%%%%%%%%%%%%%%%%%%%%%%%%%%%%%%%%%%%%%%%%%%%%%%%%%%%%%%%%%%%%%%%%%%%%%%%%%%%%%%%%%%%%%%%%%%%%%%%%%%%%%%%%%%%%%%%%%%%%%%%%%%%%
\section{Single impulsive systems}
We consider single impulsive systems without time-delays of the form
\begin{align}\label{def:imp_sys_single}
\begin{aligned}
\dot{x}(t)&=f(x(t),u(t)),\ t\neq t_k,\ k\in\N,\\
	x(t)&=g(x^-(t),u^-(t)),\ t=t_k,\ k\in\N,
\end{aligned}
\end{align}
where $t\in\R_+$, $x\in\R^{N}$ is absolutely continuous between impulses, $u\in\R^{M}$ is a locally bounded, Lebesgue-measurable input and $\{t_1,t_2,t_3,\ldots\}$ is a strictly increasing sequence of impulse times in $(t_0,\infty)$ for some initial time $t_0< t_1$.
%, where without loss of generality we set $t_0=0$. 
The set of impulse times is assumed to be either finite or infinite and unbounded and impulse times $t_k$ have no finite accumulation point. Given a sequence $\{t_k\}$ and a pair of times $s,\ t$ satisfying $t_0\leq s<t$, $N(t,s)$ denotes the number of impulse times $t_k$ in the semi-open interval $(s,t]$.

Furthermore, $f:\R^N\times\R^{M}\rightarrow\R^{N}$, $g:\R^N\times\R^{M}\rightarrow \R^N$, where we assume that $f$ is locally Lipschitz. All signals ($x$ and inputs $u$) are assumed to be right-continuous and to have left limits at all times and we denote $x^{-}:=\lim_{s\nearrow t}x(s)$, $u^{-}(t):=\lim_{s\nearrow t}u(s)$.

%\makebox[-20pt][r]{\textbf{Kommentar!!!}}
%\fbox{Hier Änderungen: supremum norm von $u$ rausgenommen. 
%}

%The supremum norm of an input $u$ on the interval $[t_0,t]$ is defined by
%$$\|u\|_{[t_0,t]}:=\max\left\{\mathrel{\mathop{\text{ess sup}}\limits_{ s\in[t_0,t]}}|u(s)|,\sup\limits_{t_k\in[t_0,t]}|u(t_k)|\right\},$$
%where by $t_k\in[t_0,t]$ we consider impulse times $t_k$ between $t_0$ and $t$.
We are interested in the stability of systems of the form \eqref{def:imp_sys_single}, where we use the following stability property, introduced in \cite{Son89} and adapted to impulsive systems in \cite{HLT08} as follows:
\begin{definition}
Assume that a sequence $\{t_{k}\}$ is given. We call system \eqref{def:imp_sys_single} {\it input-to-state stable (ISS)} if there exist functions $\beta\in\KL$, $\gamma\in\K_{\infty}$, such that for every initial condition $x(t_0)$ and every input $u$ the corresponding solution to \eqref{def:imp_sys_single} exists globally and satisfies
\begin{equation}\label{eq:iss_ws}
|x(t)|\leq\max\{\beta(|x(t_0)|,t-t_0),\gamma(\|u\|_{[t_0,t]})\},\quad\forall t\geq t_0.
\end{equation}
The impulsive system \eqref{def:imp_sys_single} is {\it uniformly ISS} over a given class $\Simp$ of admissible sequences of impulse times if (\ref{eq:iss_ws}) holds for every sequence in $\Simp$, with functions $\beta$ and $\gamma$ that are independent of the choice of the sequence.
\end{definition}
For the stability analysis of impulsive systems we use exponential ISS-Lyapunov functions, see \cite{HLT08}. Here, we assume that these functions are locally Lipschitz continuous, which are differentiable for almost all (f.a.a.) $x$ by Rademacher's Theorem (see e.g., \cite{Eva98} Theorem 5.8.6). For the stability analysis of interconnected systems it is sufficient to consider locally Lipschitz continuous functions instead of smooth functions and they were also used for example in \cite{DRW10}.
\begin{definition}\label{def:iss_lyap_ws}
We say that a function $V:\R^N\rightarrow\R_+$ is an {\it exponential ISS-Lyapunov function} for \eqref{def:imp_sys_single} with {\it rate coefficients} $c, d\in \R$ if $V$ is locally Lipschitz, positive definite, radially unbounded, and whenever $V(x)\geq\gamma(|u|)$ holds it follows
\begin{align}
\nabla V(x)\cdot f(x,u)\leq -cV(x)\ \text{f.a.a. }x,\text{ all }u\text{ and}\label{eq:iss_lyap1_ws}\\
V(g(x,u))\leq e^{-d}V(x)\ \forall x, u,\label{eq:iss_lyap2_ws}
\end{align}
where $\gamma$ is some function from $\K_{\infty}$.
\end{definition}
%Roughly speaking, condition \eqref{eq:iss_lyap1_ws} states, that if $c$ is positive then the function $V$ decreases along the trajectory. On the other hand, if $c<0$ then the function $V$ can increase along the trajectory. Condition \eqref{eq:iss_lyap2_ws} states, that if $d$ is positive then the jump (impulse) decreases the magnitude of $V$. On the other hand, if $d<0$ then the jump (impulse) can increase the magnitude of $V$.

Without loss of generality we use the same function $\gamma$ in \eqref{eq:iss_lyap1_ws} and \eqref{eq:iss_lyap2_ws}. Choosing ${\gamma}_c\in\K_{\infty}$ in \eqref{eq:iss_lyap1_ws} and ${\gamma}_t\in\K_{\infty}$ in \eqref{eq:iss_lyap2_ws} and taking the maximum of these two functions, we get $\gamma$.

Note that conditions \eqref{eq:iss_lyap1_ws} and \eqref{eq:iss_lyap2_ws} are in implication form. This is equivalent to the usage of the dissipative form in \cite{HLT08}, which was proved in \cite{CaT09}, Proposition~2.6., where the coefficients $c$, $d$ are different in general.

In \cite{HLT08} the following theorem was proved which establishes stability of a single impulsive system.
\begin{theorem}[Lyapunov-type theorem]\label{thm:stab_1sys}
Let $V$ be an exponential ISS-Lyapunov function for (\ref{def:imp_sys_single}) with rate coefficients $c, d\in \R$ with $d\neq 0$. For arbitrary constants $\mu, \lambda>0$, let $\Simp[\mu,\lambda]$ denote the class of impulse time sequences $\{t_k\}$ satisfying 
\begin{align}\label{eq:cond_int}
-dN(t,s)-(c-\lambda)(t-s)\leq\mu,\ \forall t\geq s \geq t_0.
\end{align}
Then the system (\ref{def:imp_sys_single}) is uniformly ISS over $\Simp[\mu,\lambda]$.
\end{theorem}

\begin{remark}
Note that in \cite{HLT08} this theorem was proved, using exponential ISS-Lyapunov functions in dissipative form. However, the same statement holds also if the exponential ISS-Lyapunov function is given in the implication form \eqref{eq:iss_lyap1_ws} and \eqref{eq:iss_lyap2_ws}. The proof is similar to the proof of Theorem~\ref{th:SingleIMPTDSISSLyaKra}.
%By Proposition~2.6 in \cite{CaT09} the conditions in dissipative form are equivalent to the conditions in implication form, used in Definition \ref{def:iss_lyap_ws}, but the coefficients $c$, $d$ may be different. 	
\end{remark}

If $d=0$, then the jumps do not destabilize the system, and the whole system will be ISS, if the corresponding continuous dynamics is ISS.
This case was investigated in more detail in \cite{HLT08}, Section~6.
 
Note, that condition \eqref{eq:cond_int} guarantees stability of the impulsive system even if the continuous or discontinuous behavior is unstable. For example, if the continuous behavior is unstable, which means $c<0$, then this condition assumes that the discontinuous behavior has to stabilize the system $(d>0)$ and the jumps have to occur often enough. Conversely, if the discontinuous behavior is unstable $(d<0)$ and the continuous behavior is stable $(c>0)$ then the jumps have to occur rarely, which stabilizes the system.

%%%%%%%%%%%%%%%%%%%%%%%%%%%%%%%%%%%%%%%%%%%%%%%%%%%%%%%%%%%%%%%%%%%%%%%%%%%%%%%%%%%%%%%%%%%%%%%%%%%%%%%%%%%%%%%%%%%%%%%%%%%%%%%%%%%%%%%%%%%%%%%%%%%%%%%%%%%%%%%%%%%
\subsection{Systems with time-delays}

We consider single impulsive system with time-delays of the form
\begin{align}\label{eq:SingleIMPTDS}
\begin{aligned}
\dot{x}(t)=&\ f(x^t,u(t)),\ t\neq t_k,\ k\in\N,\\
	x(t)=&\ g((x^t)^-,u^-(t)),\ t=t_k,\ k\in\N,
\end{aligned}
\end{align}
where we make the same assumptions as in the delay-free case, where $f:PC\left(\left[-\theta,0\right];\R^N\right)\times\R^{M}\rightarrow\R^{N}$ is locally Lipschitz and $g:PC\left(\left[-\theta,0\right];\R^N\right)\times\R^{M}\rightarrow \R^N$. We denote $(x^t)^{-}:=\lim_{s\nearrow t}x_s$.

We assume that the regularity conditions (see e.g., \cite{BaX00}) for the existence and uniqueness of a solution of system \eqref{eq:SingleIMPTDS} are satisfied. We denote the solution of \eqref{eq:SingleIMPTDS} corresponding to the given input $u$ by $x(t,t_0,\xi,u)$ or $x(t)$ for short, for any $\xi\in PC([-\theta,0],\R^{N})$ that exists in a maximal interval $[-\theta,b)$, $0<b\leq+\infty$, satisfying the initial condition $x^{t_0}=\xi$.

In the following subsections we present tools, namely Lyapunov-Krasovskii functionals and Lyapunov-Razumikhin functions, to check, if a time-delay system has the ISS property.

\subsubsection{The Lyapunov-Krasovskii methodology}\label{subsec:LKfunctionalsingle}
In this subsection we adapt the Lyapunov-Krasovskii methodology, introduced in \cite{PeJ06}, to impulsive time-delay systems. As one of the results of this paper we prove that from the existence of an exponential ISS-Lyapunov-Krasovskii functional the ISS property follows, provided that the dwell-time condition \eqref{eq:cond_int} is satisfied.

We consider another type of impulsive systems with time-delays of the form
\begin{align}\label{eq:SingleIMPTDSstatejump}
\begin{aligned}
\dot{x}(t)=&\ f(x^t,u(t)),\ t\neq t_k,\ k\in\N,\\
	x^t=&\ g((x^t)^-,u^-(t)),\ t=t_k,\ k\in\N,
\end{aligned}
\end{align}
where we make the same assumptions as before and the functional $g$ is now a map from $PC\left(\left[-\theta,0\right];\R^N\right)\times\R^{M}$ into $PC\left(\left[-\theta,0\right];\R^N\right)$.

According to \cite{HiP05}, Section 2, the initial state and the input together determine the evolution of the system according to the right-hand side of the differential equation. Therefore, for time-delay systems we denote the state by the function $x^t\in PC([-\theta,0],\R^N)$ and we change the discontinuous behavior in \eqref{eq:SingleIMPTDSstatejump}: In contrary to the system \eqref{eq:SingleIMPTDS} at an impulse time $t_k$ not only the point $x(t_k)$ ``jumps'', but all the states $x^t$ in the interval $(t_k-\theta,t_k]$. Due to this change the Lyapunov-Razumikhin approach cannot be applied. In this case we propose to use Lyapunov-Krasovskii functionals for the stability analysis of systems of the form \eqref{eq:SingleIMPTDSstatejump}.

Another approach using Lyapunov functionals can be found in \cite{SLL99}. There, Lyapunov functionals for systems of the form \eqref{eq:SingleIMPTDS} with zero input are used for stabilization of impulsive systems, where the definition of such a functional is different to the approach presented here according to impulse times.

The ISS property is redefined with respect to time-delays, see \cite{ChZ09}:
\begin{definition}\label{def:ISS}
Suppose that a sequence $\{t_k\}$ is given. We call system \eqref{eq:SingleIMPTDS} (or \eqref{eq:SingleIMPTDSstatejump}) {\it input-to-state stable (ISS)} if there exist functions $\beta\in\KL$, $\gamma_u\in\K_{\infty}$, such that for every initial condition $\xi \in PC([-\theta,0],\R^N)$ and every input $u$ the corresponding solution to \eqref{eq:SingleIMPTDS} (or \eqref{eq:SingleIMPTDSstatejump}) exists globally and satisfies
\begin{align}\label{eq:iss_imp_delay}
|x(t)|\leq\max\{\beta(\|\xi\|_{[-\theta,0]},t-t_0),\gamma_u(\|u\|_{[t_0,t]})\},\ \forall t\geq t_0.
\end{align}
The impulsive system \eqref{eq:SingleIMPTDS} (or \eqref{eq:SingleIMPTDSstatejump}) is {\it uniformly ISS} over a given class $\Simp$ of admissible sequences of impulse times if \eqref{eq:iss_imp_delay} holds for every sequence in $\Simp$, with functions $\beta$ and $\gamma_u$, which are independent on the choice of the sequence.
\end{definition}

Given a locally Lipschitz continuous functional $V:PC\left(\left[-\theta,0\right];\R^N\right)\rightarrow\R_+$, the upper right-hand derivative $D^+V$ of the functional $V$ along the solution $x(t,t_0,\xi,u)$ is defined according to \cite{HaV93}, Chapter 5.2:
\begin{align}\label{eq:AblLK}
D^+V\left(\phi,u\right):=\limsup_{h\rightarrow 0^+}\tfrac{1}{h}\left(V\left(x^{t+h}\right)-V\left(\phi\right)\right),
\end{align}
where $x^{t+h}\in PC\left(\left[-\theta,0\right];\R^N\right)$ is generated by the solution $x(t,t_0,\phi,u)$ of the system \eqref{eq:SingleIMPTDSstatejump} with $x^{t_0}:=\phi \in PC\left(\left[-\theta,0\right];\R^N\right)$.

With the symbol $|\cdot|_a$ we indicate any norm in $PC\left(\left[-\theta,0\right];\R^N\right)$ such that for some positive reals $b,\tilde{c}$ the following inequalities hold
\begin{align*}
b|\phi(0)|\leq|\phi|_a\leq \tilde{c}\|\phi\|_{[-\theta,\infty)},\ \forall \phi\in PC\left(\left[-\theta,0\right];\R^N\right).
\end{align*}

\begin{definition}\label{def:ISS_LK_functional}
A functional $V:PC\left(\left[-\theta,0\right];\R^N\right)$ $\rightarrow\R_+$ is called an \textsl{exponential ISS-Lyapunov-Krasovskii functional} with rate coefficients $c,d\in\R$ for system \eqref{eq:SingleIMPTDSstatejump}, if $V$ is locally Lipschitz continuous, there exist $\psi_1,\psi_2\in{\cK}_{\infty}$ such that
\begin{align}\label{eq:BedLKfunctional1}
	\psi_1(|\phi(0)|)\leq V(\phi) \leq \psi_2(|\phi|_a),\ \forall \phi\in PC\left(\left[-\theta,0\right];\R^N\right)
\end{align}
and there exists a function $\gamma\in{\cK}$ such that whenever $V(\phi)\geq \gamma(|u|)$ holds it follows
\begin{align}
	\D^+V(\phi,u)\leq-cV(\phi)\ \text{and}\label{eq:CondLK2}\\
	 V(g(\phi,u))\leq e^{-d}V(\phi),\label{eq:CondLK3}
%V(\phi)\geq \gamma(|u(t)|)&\Rightarrow \D^+V(\phi,u)\leq-cV(\phi),\label{eq:CondLK2}\\
%V(\phi)\geq \gamma(|u(t)|)&\Rightarrow V(g(\phi,u))\leq e^{-d}V(\phi),\label{eq:CondLK3}
\end{align}
for all $\phi\in PC\left(\left[-\theta,0\right];\R^N\right)$ and $u\in\R^M$.
\end{definition}
%The rate coefficients $c, d$ are not required to be non-negative and therefore $V$ may not decrease even if $u=0$.

%A different but equivalent formulation of an exponential ISS-Lyapunov-Krasovskii functional is the following:
%\begin{proposition}\label{prop:equivlance_LK}
%Assume that $g$ is continuous with $g(0,0)=0$. $V$ is exponential ISS-Lyapunov-Krasovskii functional for \eqref{eq:SingleIMPTDSstatejump} if and only if $V$ satisfies the requirements of Definition~\ref{def:ISS_LK_functional}, where \eqref{eq:CondLK3} is replaced by
%\begin{align*}
%	V(g(\phi,u))\leq \max\{e^{-d}V(\phi),\tilde{\gamma}(|u|)\}\ \forall \phi, u,%\label{eq:CondLK3NEU}
%\end{align*}
%$\tilde{\gamma}\in{\cK}_{\infty}$.
%\end{proposition}
%
%\textit{\textbf{Proof.}}
%The proof follows the same steps as in the proof of Proposition~\ref{prop:equivlance} with according changes to functionals and time-delays.
%\hfill $\Box$

Now, we present as a result a counterpart of Theorem~1 in \cite{HLT08} and Theorems~1 and 2 in \cite{ChZ09} for impulsive systems with time-delays using the Lyapunov-Krasovskii approach:

\begin{theorem}[Lyapunov-Krasovskii-type theorem]\label{th:SingleIMPTDSISSLyaKra}
Let $V$ be an exponential ISS-Lyapunov-Krasovskii functional for system \eqref{eq:SingleIMPTDSstatejump} with $c,d\in\R,\ d\neq0$. For arbitrary constants $\mu,\ \lambda\in\R_+$, let $\Simp [\mu,\lambda]$ denote the class of impulse time sequences $\{t_k\}$ satisfying the dwell-time condition \eqref{eq:cond_int}. Then the system \eqref{eq:SingleIMPTDSstatejump} is uniformly ISS over $\Simp [\mu,\lambda]$.
\end{theorem}

%\textit{\textbf{Proof.}} 
%The proof is similar to the proof of Theorem~1 in \cite{HLT08} with corresponding changes due to time-delays and the usage of a functional.
%\hfill $\Box$

\textit{\textbf{Proof.}}
From \eqref{eq:CondLK2} we have for any two consecutive impulses $t_{k-1},t_k$, $\forall t\in(t_{k-1},t_k)$
\begin{align}\label{eq:CondLK2neu}
	V(x^t)\geq \gamma(|u(t)|)\Rightarrow \D^+V(x^t,u(t))\leq-cV(x^t)
\end{align}
and similarly with \eqref{eq:CondLK3} for every impulse time $t_k$
\begin{align}\label{eq:CondLK3neu}
	V(x^{t_k})\geq \gamma(|u^-(t_k)|)\Rightarrow V(g((x^{t_k})^-,u^-(t_k)))\leq e^{-d}V(x^{t_k}).
\end{align}
Because of the right-continuity of $x$ and $u$ there exists a sequence of times $t_0:= \tilde{t}_0<\bar{t}_1<\tilde{t}_1<\bar{t}_2<\tilde{t}_2<\ldots$ such that we have
\begin{align}
	V(x^t)\geq\gamma(\|u\|_{[t_0,t]}),\ \forall t\in[\tilde{t}_i,\bar{t}_{i+1}),\ i=0,1,\ldots,\label{eq:Fall1KF}\\
	V(x^t)\leq\gamma(\|u\|_{[t_0,t]}),\ \forall t\in[\bar{t}_i,\tilde{t}_{i}),\ i=1,2,\ldots,\label{eq:Fall2KF}
\end{align}
where this sequence breaks the interval $[t_0,\infty)$ into a disjoint union of subintervals. Suppose $t_0<\bar{t}_1$, so that $[t_0,\bar{t}_1)$ is nonempty. Otherwise skip forward to the line below \eqref{eq:Bedcont}.
Between any two consecutive impulses $t_{k-1},t_k\in(t_0,\bar{t}_1]$ with \eqref{eq:Fall1KF} and \eqref{eq:CondLK2neu} we have $\D^+V(x^t,u(t))\leq-cV(x^t),\ \forall t\in(t_{k-1},t_k)$ and therefore
\begin{align*}
	V((x^{t_k})^-)\leq e^{-c(t_k-t_{k-1})}V(x^{t_{k-1}}).
\end{align*}
From \eqref{eq:CondLK3neu} and \eqref{eq:Fall1KF} we have $V(x^{t_k})\leq e^{-d}V((x^{t_k})^-)$. Combining this it follows
\begin{align*}
	V(x^{t_k})\leq e^{-d}e^{-c(t_k-t_{k-1})}V(x^{t_{k-1}})
\end{align*}
and by the iteration over the $N(t,t_0)$ impulses on $(t_0,t]$ we obtain the bound
\begin{align*}
	V(x^t)\leq e^{-dN(t,t_0)-c(t-t_0)}V(\xi),\ \forall t\in(t_0,\bar{t}_1].
\end{align*}
Using the dwell-time condition \eqref{eq:cond_int} we get
\begin{align}\label{eq:Bedcont}
	V(x^t)\leq e^{\mu-\lambda(t-t_0)}V(\xi),\ \forall t\in(t_0,\bar{t}_1].
\end{align}
Now, on any subinterval of the form $[\bar{t}_i,\tilde{t}_i)$ we already have \eqref{eq:Fall2KF} as a bound. If $\tilde{t}_i$ is not an impulse time, then \eqref{eq:Fall2KF} is a bound for $t=\tilde{t}_i$. If $\tilde{t}_i$ is an impulse time, then we have
\begin{align*}
V(x^{\tilde{t}_i})\leq e^{-d}\gamma(\|u\|_{[t_0,\tilde{t}_i]})	
\end{align*}
and in either case
\begin{align}\label{eq:Bedsprung}
	V(x^t)\leq e^{|d|}\gamma(\|u\|_{[t_0,t]}),\ \forall t\in[\bar{t}_i,\tilde{t}_i],\ i\geq 1,
\end{align}
where this bound holds for all $t\geq \bar{t}_i$, if $\tilde{t}_i=\infty$. Now consider any subinterval of the form $[\tilde{t}_i,\bar{t}_{i+1}),\ i\geq 1$. Repeating the argument used to establish \eqref{eq:Bedcont} with $\tilde{t}_i$ in place of $t_0$ and using \eqref{eq:Bedsprung} with $t=\tilde{t}_i$, we get
\begin{align*}
	V(x^t)\leq e^{\mu-\lambda(t-\tilde{t}_i)}V(x^{\tilde{t}_i})\leq e^{\mu+|d|}\gamma(\|u\|_{[t_0,\tilde{t}_i]})
\end{align*}
$\forall t\in(\tilde{t}_i,\bar{t}_{i+1}],\ i\geq 1$. Combining this with \eqref{eq:Bedcont} and \eqref{eq:Bedsprung} we obtain the global bound $\forall t\geq t_0$
\begin{align*}
	V(x^t)\leq\max\{e^{\mu-\lambda(t-t_0)}V(\xi),e^{\mu+|d|}\gamma(\|u\|_{[t_0,t]})\}.
\end{align*}
The uniformly ISS property follows then from \eqref{eq:BedLKfunctional1} by definition of $\beta(\xi,t-t_0):=\psi_1^{-1}(e^{\mu-\lambda(t-t_0)}\psi_2(\tilde{c}\|\xi\|_{[-\theta,0]}))$ and $\gamma_u(r):=\psi_1^{-1}(e^{\mu+|d|}\gamma(r))$, where the global existence of solutions follows from the boundedness of $x$. Note that $\beta$ and $\gamma_u$ do not depend on the particular choice of the time sequence and therefore uniformity is clear.\hfill $\Box$
 %%%%%%%%%%%%%%%%%%%%%%%%%%%%%%%%%%%%%%%%%%%%%%%%%%%%%%%%%%%%%%%%%%%%%%%%%%%%%%%%%%%%%%%%%%%%%%%%%%%%%%%%%%%%%%%%%%%%%%%%%%%%%%%%%%%%%%%%%%%%%%%%%%%%%%%%%%%%%%%%%%%
\subsubsection{The Lyapunov-Razumikhin methodology}\label{subsec:LRfunctionsingle}
To study the ISS property of impulsive systems with time-delays of the form \eqref{eq:SingleIMPTDS} one can use ISS-Lyapunov-Razumikhin functions.

\begin{definition}\label{DefLRfunction}
A function $V:\R^N\rightarrow\R_+$ is called an \textsl{exponential ISS-Lyapunov-Razumikhin function} for system \eqref{eq:SingleIMPTDS} with rate coefficients $c,d\in\R$, if $V$ is locally Lipschitz continuous and there exist functions $\psi_1,\psi_2$, $\gamma_t,\ \gamma_u\in{\cK}_{\infty}$ such that 
\begin{align*}
	\psi_1(|\phi(0)|)\leq V(\phi(0)) \leq \psi_2(|\phi(0)|)
\end{align*}
and whenever $V(\phi(0))\geq \max\{\gamma_t(\|V^t(\phi)\|),\gamma_u(|u|)\}$ holds it follows
\begin{align}
				 \D^+V(\phi(0))\leq-cV(\phi(0))\ \text{and}\label{CondLR2}\\
	 V(g(\phi,u))\leq e^{-d}V(\phi(0)),\label{CondLR3}
\end{align}
for all $\phi\in PC\left(\left[-\theta,0\right];\R^N\right)$ and $u\in\R^M$, where $V^t:PC\left(\left[-\theta,0\right];\R^N\right) \to PC\left(\left[-\theta,0\right];\R^N\right)$ is defined by $V^t(x^t)(\tau):=V(x(t+\tau)),\ \tau\in\left[-\theta,0\right]$.
\end{definition}

For the main result of this section, we need the following:

\begin{definition}
\label{def:GS}
Assume that a sequence $\{t_{k}\}$ is given. We call system 
\eqref{eq:SingleIMPTDS} (or \eqref{eq:SingleIMPTDSstatejump})
{\it globally stable (GS)} if there exist functions $\varphi, \gamma\in\K_{\infty}$,
such that for every initial condition $\xi \in PC([-\theta,0],\R^N)$ and every input $u$ it holds
\begin{align}
\label{eq:gs_imp_delay}
|x(t)|\leq\max\{\varphi(\|\xi\|_{[-\theta,0]}),\gamma(\|u\|_{[t_0,t]})\},\ \forall t\geq t_0.
\end{align}
The impulsive system \eqref{eq:SingleIMPTDS} (or \eqref{eq:SingleIMPTDSstatejump}) is {\it uniformly GS} over a given class $\Simp$ of admissible sequences of impulse times if \eqref{eq:gs_imp_delay} holds for every sequence in $\Simp$, with functions $\varphi$ and $\gamma_u$, which are independent on the choice of the sequence.
\end{definition}

To prove the Razumikhin-type theorem for impulsive time-delay systems, we need the following characterization of the uniform ISS property:

\begin{lemma}
\label{UISS_Characterisation}
The system \eqref{eq:SingleIMPTDS} (or \eqref{eq:SingleIMPTDSstatejump}) is uniformly ISS over $\Simp$ if and only if it is
\begin{itemize}
	\item uniformly globally stable over $\Simp$ and
	\item $\exists \gamma \in \K$ such that for each $\epsilon>0$, $\eta_x\in\R_+,\ \eta_u\in\R_+$ there exists $T\geq0$ (which does not depend on an impulse time sequence from $\Simp$) such that $\|\xi\|_{[-\theta,0]}\leq\eta_x$ and $\|u\|_{\infty}\leq\eta_u$ imply $|x(t)|\leq\max\{\epsilon,\gamma(\|u\|_{[t_0,t]})\},\ \forall t \geq T + t_0$.
\end{itemize}
\end{lemma}

\textit{\textbf{Proof.}}
We start with necessity. Let \eqref{eq:SingleIMPTDS} (or \eqref{eq:SingleIMPTDSstatejump}) be uniformly ISS over $\Simp$. Then it is uniformly GS over $\Simp$ with $\varphi(\cdot):=\beta(\cdot,0),\ \gamma_u\equiv\gamma$. 

Take arbitrary $\eps>0$, $\eta_x\in\R_+$. 
For all $\|\xi\|_{[-\theta,0]}\leq\eta_x$, all $u$ and all impulse time sequences from $\Simp$ it holds
\begin{align*}
|x(t)|\leq\max\{\beta(\eta_x,t-t_0),\gamma_u(\|u\|_{[t_0,t]})\},\ \forall t\geq t_0.
\end{align*}
If $\eps > \beta(\eta_x,0)$, then we choose $T$ in Lemma~\ref{UISS_Characterisation} as $T:=0$. 
Otherwise take $T$ as solution (which is unique) of the equation $\beta(\eta_x,T-t_0)=\eps$. Clearly, $T$ does not depend on the choice of the sequence from $\Simp$. Thus, the second property in the statement of the lemma is verified with $\gamma_u\equiv\gamma$.

Now let us prove sufficiency. Without loss of generality we take $\kappa:=\eta_x=\eta_u$ and fix it. 
From uniform global stability it follows that for all $\|\xi\|_{[-\theta,0]}\leq\kappa$, for all $u \in L_{\infty}(\R_+,\R^m)$ and all impulse time sequences from $\Simp$ it holds
\begin{align*}
|x(t)|\leq\max\{\varphi(\kappa),\gamma(\|u\|_{[t_0,t]})\},\ \forall t\geq t_0.
\end{align*}
Define $\eps_n:=\tfrac{1}{2^n}\varphi(\kappa)$. The second assumption of the lemma implies the existence of a sequence of times $T_n:=T(\eps_n,\kappa)$, which without loss of generality we assume to be strictly increasing such that for all
$\xi: \|\xi\|_{[-\theta,0]}\leq \kappa,$ for all $u: \|u\|_{\infty}\leq \kappa$ and for all $t\geq T_n+t_0$ it holds
\begin{align*}
|x(t)|\leq \max\left\{\eps_n,\gamma(\|u\|_{[t_0,t]})\right\}.
\end{align*} 
Define $\omega(\kappa,T_n):=\eps_{n-1}$, $n\in\N,n\neq0$. Extend this function for $t\in\R_+\backslash\{T_n,n\in\N\}$ such that $\omega(\kappa,\cdot)\in{\cL}$. For all $t\in(T_n+t_0,T_{n+1}+t_0)$ it holds
\begin{align*}
	|x(t)|\leq\max\left\{\eps_n,\gamma(\|u\|_{[t_0,t]})\right\}\leq\max\left\{\omega(\kappa,t-t_0),\gamma(\|u\|_{[t_0,t]})\right\}.
\end{align*} 
Doing this for all $\kappa\in\R_+$ we obtain the function $\omega$, defined on $\R_+\times\R_+$. Now define $\beta(r,t)=\sup_{0\leq s \leq r}\omega(s,t)\geq\omega(r,t)$. It follows that $\beta$ is continuous, $\beta(\cdot,t)\in{\cK}$ (if it is not true, we can always find $\tilde{\beta}(\cdot,t)\in\K$ without losing continuity in the second argument such that $\beta(r,t)\leq\tilde{\beta}(r,t)$ for all $r>0$) and $\beta(r,\cdot)\in{\cL}$, since $\omega(r,\cdot)\in{\cL}$. Thus, $\beta\in{\cK\cL}$ and 
%by definition of $\gamma(r):=\psi_1^{-1}(e^{\mu}\gamma_u(r)),\ r\geq0$ 
we obtain
\begin{align*}
	|x(t)| & \leq \max\left\{\beta(\max\{ \|\xi\|_{[-\theta,0]},\|u\|_{[t_0,t]} \},t-t_0),\gamma(\|u\|_{[t_0,t]})\right\} \\
	 & = \max\left\{ \beta(\|\xi\|_{[-\theta,0]},t-t_0),\beta(\|u\|_{[t_0,t]},t-t_0),\gamma(\|u\|_{[t_0,t]})\right\} \\
	& \leq \max\left\{ \beta(\|\xi\|_{[-\theta,0]},t-t_0),\gamma_u(\|u\|_{[t_0,t]})\right\},	
\end{align*}
where $\gamma_u(r):=\max\{\beta(r,0), \gamma(r)\}$. This proves uniform ISS over $\Simp$. \hfill $\Box$

%The following theorem is similar to Theorem~\ref{th:SingleIMPTDSISSLyaKra} in this paper and Theorem~1 in \cite{HLT08}. In \cite{ChZ09} non-autonomous systems are investigated and a different approach for the characterization of the condition on the time intervals of the impulses is proposed. There, a so-called \textsl{fixed} dwell-time condition and a Lyapunov-Razumikhin approach are used. Here, we consider the \textsl{average} dwell-time condition, which was also used in \cite{HLT08} and introduced in \cite{HeM99} for switched systems. The average dwell-time condition considers the average of impulses over an interval, whereas the fixed dwell-time condition considers the (minimal or maximal) interval between two impulses, which lead to a more conservative condition in contrast to the average dwell-time condition.

\begin{remark}
Note that in \cite{Tee98} for time-delay systems without impulse times a different definition of the ISS property has been used. However, with the help of Lemma~\ref{UISS_Characterisation} one can show that from the Assumptions of Theorem~1 in \cite{Tee98} it follows not only ISS in the sense of \cite{Tee98}, but also ISS in the sense of Definition~\ref{def:ISS}. Until now, it is an open problem whether the definition of ISS in \cite{Tee98} is equivalent to the ISS property in Definition~\ref{def:ISS} for time-delay systems without impulse times, see \cite{TWJ09,TiW10}.  
\end{remark}

We need one more technical result for the main result of this section. For a given sequence of impulse times, we denote the number of jumps within the time-span $[s,t]$ by $N^*(t,s)$. The set of impulse time sequences, for which \eqref{eq:cond_int}
% \eqref{Classical_Dwell-Time-Cond} 
holds with $N^*(t,s)$ instead of $N(t,s)$, is denoted by $\Simp^*[\mu,\lambda]$. 
The xext lemma shows that $\Simp^*[\mu,\lambda]$ is equal to $\Simp[\mu,\lambda]$.
\begin{lemma}
\label{Dwell-Time_Equiv}
Let $c,d \in \R$, $d \neq 0$ be given. Then, $\Simp[\mu,\lambda] = \Simp^*[\mu,\lambda]$ for all $\mu,\lambda>0$.
\end{lemma}

\textit{\textbf{Proof.}}
Firstly, consider the case, when $d<0$.
Let $T \in \Simp^*[\mu,\lambda]$. Since $N(t,s) \leq N^*(t,s) \leq \frac{1}{-d}( (c - \lambda)(t-s) + \mu)$, it is clear that $\Simp[\mu,\lambda] \supset \Simp^*[\mu,\lambda]$.

Let $T \in \Simp[\mu,\lambda]$. Then, since $[s,t] \subset (s-\eps,t]$ for arbitrary $\eps>0$ it holds
\[
N^*(t,s) \leq N(t,s-\eps) \leq \frac{1}{-d}( (c - \lambda)(t-(s-\eps)) + \mu).
\]
Tending $\eps \to 0$, we obtain that $\Simp[\mu,\lambda] \subset \Simp^*[\mu,\lambda]$.

Now let $d>0$. 
Take $T \in \Simp[\mu,\lambda]$. Then, $N^*(t,s) \geq  N(t,s)$ implies $\Simp[\mu,\lambda] \subset \Simp^*[\mu,\lambda]$. 
On the other side, since for arbitrary $\eps>0$ it holds
\[
N(t,s) \geq N^*(t,s+\eps) \geq \frac{1}{-d}( (c - \lambda)(t-(s+\eps)) + \mu),
\]
for $T \in \Simp^*[\mu,\lambda]$, tending $\eps \to 0$ we have $\Simp[\mu,\lambda] \supset \Simp^*[\mu,\lambda]$. \hfill $\Box$
%\end{proof}

Now, we prove the main result of this section.

\begin{theorem}[Lyapunov-Razumikhin-type theorem]\label{PropISSRazumikhin}
Let $V$ be an exponential ISS-Lyapunov-Razumikhin function for system \eqref{eq:SingleIMPTDS} with $c,d\in\R,\ d\neq0$. For arbitrary constants $\mu,\ \lambda >0$, let $\Simp [\mu,\lambda]$ denote the class of impulse time sequences $\{t_k\}$ satisfying the dwell-time condition \eqref{eq:cond_int}. If $\gamma_t$ satisfies $\gamma_t(r)<e^{-\mu}r,\ r>0$, then the system \eqref{eq:SingleIMPTDS} is uniformly ISS over $\Simp [\mu,\lambda]$.%$\text{ISS}^{\text{T}}$
\end{theorem}

\textit{\textbf{Proof.}}
From \eqref{CondLR2} we have for any two consecutive impulses $t_{k-1},t_k$, $\forall t\in(t_{k-1},t_k)$
\begin{align}\label{eq:CondLR2neu}
	V(x(t))\geq \max\{\gamma_t(\|V^t(x^t)\|),\gamma_u(|u(t)|)\} \Rightarrow \D^+V(x(t))\leq-cV(x(t))
\end{align}
and similarly with \eqref{CondLR3} for every impulse time $t_k$
\begin{eqnarray}
\label{eq:CondLR3neu}
	V(x(t_k))\geq &\max\{\gamma_t(\|V^t((x^{t_k})^-)\|),\gamma_u(|u^-(t_k)|)\} \\
	&\Rightarrow V(g((x^{t_k})^-,u^-(t_k)))\leq e^{-d}V(x(t_k)). \notag
\end{eqnarray}
Because of the right-continuity of $x$ and $u$ there exists a sequence of times $t_0:= \tilde{t}_0<\bar{t}_1<\tilde{t}_1<\bar{t}_2<\tilde{t}_2<\ldots$ such that for $i=0,1,\ldots$ we have
\begin{align}
\label{eq:Fall1}
	V(x(t))\geq\max\{\gamma_t(\sup_{r\in\left[t_0,t\right]}\left\|V^t(x^r)\right\|),\gamma_u(\|u\|_{[t_0,t]})\},\ \forall t\in[\tilde{t}_i,\bar{t}_{i+1})
\end{align}
and  for all $i=1,2,\ldots$ it holds
\begin{align}
\label{eq:Fall2}
V(x(t))\leq\max\{\gamma_t(\sup_{r\in\left[t_0,t\right]}\left\|V^t(x^r)\right\|),\gamma_u(\|u\|_{[t_0,t]})\},\ \forall t\in[\bar{t}_i,\tilde{t}_{i}),
\end{align}
where this sequence breaks the interval $[t_0,\infty)$ into a disjoint union of subintervals. Suppose $t_0<\bar{t}_1$, so that $[t_0,\bar{t}_1)$ is nonempty. Otherwise skip forward to the line below \eqref{eq:BedcontLR}.
Between any two consecutive impulses $t_{k-1},t_k\in[t_0,\bar{t}_1]$ from \eqref{eq:Fall1} and \eqref{eq:CondLR2neu} we have $\D^+V(x(t)))\leq-cV(x(t)),\ \forall t\in(t_{k-1},t_k)$ and therefore
\begin{align*}
	V(x^-(t_k))\leq e^{-c(t_k-t_{k-1})}V(x(t_{k-1})).
\end{align*}
From \eqref{eq:CondLR3neu} and \eqref{eq:Fall1} we have $V(x(t_k))\leq e^{-d}V(x^-(t_k))$. Combining this it follows
\begin{align*}
	V(x(t_k))\leq e^{-d-c(t_k-t_{k-1})}V(x(t_{k-1}))
\end{align*}
and by the iteration over the $N(t,t_0)$ impulses on $[t_0,t]$ we obtain the bound
\begin{align}
\label{AbschaetzungMit_N}
	V(x(t))\leq e^{-dN(t,t_0)-c(t-t_0)}V(x(t_0)),\ \forall t\in[t_0,\bar{t}_1].
\end{align}
Using the dwell-time condition \eqref{eq:cond_int} we get
\begin{align}\label{eq:BedcontLR}
	V(x(t))\leq e^{\mu-\lambda(t-t_0)}V(x(t_0)),\ \forall t\in[t_0,\bar{t}_1].
\end{align}
For any subinterval of the form $[\bar{t}_i,\tilde{t}_{i}),\ i=1,2,\ldots$ we have \eqref{eq:Fall2} as a bound for $V(x(t))$. Now consider two cases. 

Let $\tilde{t}_{i}$ be not an impulse time, then \eqref{eq:Fall2} is a bound for $t=\tilde{t}_{i}$. 
Consider the subinterval $[\tilde{t}_{i},\bar{t}_{i+1})$. Repeating the argument used to establish \eqref{eq:BedcontLR}, with $\tilde{t}_{i}$ in place of $t_0$ and using
\eqref{eq:Fall2} with $t=\tilde{t}_{i}$ we get $\forall t\in(\tilde{t}_{i},\bar{t}_{i+1}]$
\begin{align*}
	V(x(t))\leq e^{\mu-\lambda(t-\tilde{t}_{i})}V(x(\tilde{t}_{i}))\leq e^{\mu}\max\{\gamma_t (\sup_{r\in\left[t_0,\tilde{t}_{i}\right]}\left\|V^t(x^r)\right\|),\gamma_u(\|u\|_{[t_0,\tilde{t}_{i}]})\}.
\end{align*}

Now let $\tilde{t}_{i}$ be an impulse time. Then we have
\begin{align}
\label{AbschMitImpZeit}
	V(x(\tilde{t}_{i}))\leq e^{-d}\max\{\gamma_t(\sup_{r\in\left[t_0,\tilde{t}_{i}\right]}\left\|V^t(x^r)\right\|),\gamma_u(\|u\|_{[t_0,\tilde{t}_{i}]})\}.
\end{align}
Consider for one more time the subinterval of the form $[\tilde{t}_{i},\bar{t}_{i+1})$. 
According to the estimate \eqref{AbschaetzungMit_N} we obtain 
\begin{align*}
%\label{HilfAbsch_Raz}
V(x(t))\leq e^{-dN(t,\tilde{t}_{i})-c(t-\tilde{t}_{i})}V(x(\tilde{t}_{i})) 
= & e^{-d (N(t,\tilde{t}_{i})+1-1)-c(t-\tilde{t}_{i})}V(x(\tilde{t}_{i})) \\
= & e^d e^{-d N^*(t,\tilde{t}_{i})-c(t-\tilde{t}_{i})}V(x(\tilde{t}_{i})).
\end{align*}
From Lemma~\ref{Dwell-Time_Equiv} we know that $\Simp[\mu,\lambda] = \Simp^*[\mu,\lambda]$. Thus, we can continue:
\begin{align*}
%\label{HilfAbsch_Raz}
V(x(t)) \leq e^d e^{\mu-\lambda(t-\tilde{t}_{i})}V(x(\tilde{t}_{i})) 
\leq e^{\mu} \max\{\gamma_t(\sup_{r\in\left[t_0,\tilde{t}_{i}\right]}\left\|V^t(x^r)\right\|),\gamma_u(\|u\|_{[t_0,\tilde{t}_{i}]})\}.
\end{align*}

Overall, we obtain $\forall t\geq t_0$
\begin{align}\label{Absch}
	V(x(t))\leq \max\left\{e^{\mu-\lambda(t-t_0)}V(x(t_0)),e^{\mu}\gamma_t(\sup_{r\in\left[t_0,t\right]}\left\|V^t(x^r)\right\|),e^{\mu}\gamma_u(\|u\|_{[t_0,t]})\right\}
\end{align}
and it holds
\begin{align}\label{VdAbsch}
\sup_{t\geq s\geq t_0}	\left\|V^t(x^s)\right\| \leq \max\left\{\left\|V^t(x^{t_0})\right\|,\sup_{t\geq s\geq t_0}V(s)\right\}.
\end{align}
We take the supremum over $[t_0,t]$ in \eqref{Absch} and insert it into \eqref{VdAbsch}. Then, using $\gamma_t(r)<e^{-\mu}r$ and the fact that for all $a,b >0$ from $a\leq\max\{b,e^{\mu}\gamma_t(a)\}$ it follows $a\leq b$, we obtain
\begin{align*}
	\sup_{t\geq s\geq t_0}\|V^t(x^s)\|\leq\max\left\{e^{\mu}\psi_2(\|\xi\|_{[-\theta,0]}),e^{\mu}\gamma_u(\|u\|_{[t_0,t]})\right\}
\end{align*}
and therefore
\begin{align*}
	|x(t)| \leq \max\left\{\psi_1^{-1}(e^{\mu}\psi_2(\|\xi\|_{[-\theta,0]})),\psi_1^{-1}(e^{\mu}\gamma_u(\|u\|_{[t_0,t]}))\right\},\ \forall t \geq t_0,
\end{align*}
which means, that the system \eqref{eq:SingleIMPTDS} is uniformly GS over $\Simp [\mu,\lambda]$. Note that $\tilde{\varphi}( \cdot):=\psi_1^{-1}(e^{\mu}\psi_2(\cdot))$ is a ${\cK}_{\infty}$-function. Now, for given $\epsilon,\eta_x,\eta_u>0$ such that $\|\xi\|_{[-\theta,0]}\leq \eta_x,\ \|u\|_{\infty}\leq \eta_u$ let $\kappa:=\max\left\{e^{\mu}\psi_2(\eta_x),e^{\mu}\gamma_u(\eta_u)\right\}$. It holds $\sup_{t\geq s\geq t_0}\|V^t(x^s)\|\leq\kappa$. Let $\rho_2>0$ be such that $e^{-\lambda \rho_2}\kappa \leq \psi_1(\epsilon)$ and let $\rho_1>\theta$. Then, by the estimate \eqref{Absch} we have 
\begin{align*}
\sup_{t\geq s\geq t_0+\rho_1+\rho_2}\|V^t(x^s)\| &\ \leq \sup_{t\geq s \geq t_0+\rho_2}V(x(s))\\
	&\ \leq \max\left\{\psi_1(\epsilon),e^{\mu}\gamma_t(\sup_{s\in\left[t_0,t\right]}\left\|V^t(x^s)\right\|),e^{\mu}\gamma_u(\|u\|_{[t_0,t]})\right\}.
\end{align*}
In the previous inequality we put $t_0 + \rho_1+\rho_2$ instead of $t_0$ and obtain
\begin{align*}
&\ \sup_{t\geq s\geq t_0+2(\rho_1+\rho_2)}\|V^t(x^s)\|\\ 	
 \leq&\ \max\left\{\psi_1(\epsilon),e^{\mu}\gamma_t(\sup_{s\in\left[t_0+\rho_1+\rho_2,t\right]}\left\|V^t(x^s)\right\|),e^{\mu}\gamma_u(\|u\|_{[t_0 + \rho_1+\rho_2,t]})\right\} \\
 \leq&\ \max\left\{\psi_1(\epsilon),(e^{\mu}\gamma_t)^2(\sup_{s\in\left[t_0,t\right]}\left\|V^t(x^s)\right\|),e^{\mu}\gamma_u(\|u\|_{[t_0,t]})\right\}.
\end{align*}
Since $e^{\mu}\gamma_t < id$, there exists a number $\tilde{n}\in\N$, which depends on $\kappa$ and $\epsilon$ such that 
\begin{align*}
(e^{\mu}\gamma_t)^{\tilde{n}}(\kappa):=\underbrace{(e^{\mu}\gamma_t) \circ \ldots \circ (e^{\mu}\gamma_t)}_{\tilde{n} \mbox{ times}}(\kappa) \leq \max\left\{\psi_1(\epsilon),e^{\mu}\gamma_u(\|u\|_{[t_0,t]})\right\}.
\end{align*}
By induction we conclude that
\begin{align*}
	\sup_{t\geq s\geq t_0+\tilde{n}(\rho_1 + \rho_2)}\|V^t(x^s)\| \leq \max\left\{\psi_1(\epsilon),e^{\mu}\gamma_u(\|u\|_{[t_0,t]})\right\},
\end{align*}
and finally we obtain
\begin{align}\label{eq:Uatrr}
	|x(t)|\leq \max\left\{\epsilon,\psi_1^{-1}(e^{\mu}\gamma_u(\|u\|_{[t_0,t]}))\right\},\ \forall t\geq t_0+\tilde{n}(\rho_1+\rho_2).
\end{align}
Thus, the system \eqref{eq:SingleIMPTDS} satisfies the second property from the Lemma~\ref{UISS_Characterisation}, which implies that \eqref{eq:SingleIMPTDS} is uniformly ISS over $\Simp [\mu,\lambda]$.\hfill $\Box$

%Therefore, the system \eqref{eq:SingleIMPTDS} has the $\text{ISS}^{\text{T}}$ property (in the sense of Teel). The global existence of solutions follows from the boundedness of $x$. Note that $\epsilon,\gamma_u$ do not depend on the particular choice of the time sequence and therefore uniformity is clear.\hfill $\Box$

\begin{remark}
Another Razumikhin-type theorem (for non-autonomous systems) has been proposed in \cite{ChZ09}. In that paper it was used so-called \textsl{fixed} dwell-time condition to characterize the class of impulse time sequences, over which the system is uniformly ISS. In contrast to Theorems 1,2 from \cite{ChZ09}, we prove the Razumikhin-type theorem \ref{PropISSRazumikhin} over the class of sequences, which satisfy the average dwell-time condition, which is larger than the class of sequences, which satisfy the fixed dwell-time condition. However, the small-gain condition, that we have used in this paper, $\gamma_t(r)<e^{-\mu}r,\ r>0$, is stronger, than that from \cite{ChZ09}.
\end{remark}

In the next section we investigate general networks of impulsive systems in view of stability and establish a dwell-time condition for such interconnections.
%%%%%%%%%%%%%%%%%%%%%%%%%%%%%%%%%%%%%%%%%%%%%%%%%%%%%%%%%%%%%%%%%%%%%%%%%%%%%%%%%%%%%%%%%%%%%%%%%%%%%%%%%%%%%%%%%%%%%%%%%%%%%%%%%%%%%%%%%%%%%%%%%%%%%%%%%%%%%%%%%%%
\section{Large-scale networks of impulsive systems}
We consider an interconnection of $n$ impulsive subsystems with inputs of the form
%\begin{align}\label{def:imp_sys_is}
%\begin{aligned}
%\dot{x}_i(t)&=f_i(x_1(t),\ldots,x_n(t),u_i(t)),\ t\neq t_{i_{\tilde{k}}},\\
%x_i(t)&=g_i(x_1^{-}(t),\ldots,x_n^{-}(t),u_i^{-}(t)),\ t=t_{i_{\tilde{k}}},
%\end{aligned}
%\end{align}

\begin{align}\label{def:imp_sys_is}
\begin{aligned}
\dot{x}_i(t)&=f_i(x_1(t),\ldots,x_n(t),u_i(t)),\ t\neq t_k,\\
x_i(t)&=g_i(x_1^{-}(t),\ldots,x_n^{-}(t),u_i^{-}(t)),\ t=t_k,
\end{aligned}
\end{align}
$k\in\N,\ i=1,\ldots,n$, where the state $x_i(t)\in\R^{N_i}$ of the $i$th subsystem is absolutely continuous between impulses; $u_i(t)\in\R^{M_i}$ is a locally bounded, Lebesgue-measurable input and $x_j(t)\in\R^{N_j},\ j\neq i$ can be interpreted as internal inputs of the $i$th subsystem. Note that the impulse sequences for all subsystems are assumed to be equal.
%
%Given a sequence $\{t_{i_{\tilde{k}}}\}$ and a pair of times $s, t$ satisfying $t_0\leq s<t$, $N_i(t,s)$ denotes the number of impulse times $t_{i_{\tilde{k}}}$ in the semi-open interval $(s,t]$ of the $i$th subsystem.

Furthermore, $f_i:\R^{N_1}\times\ldots\times\R^{N_n}\times\R^{M_i}\rightarrow\R^{N_i}$ and $g_i:\R^{N_1}\times\ldots\times\R^{N_n}\times\R^{M_i}\rightarrow\R^{N_i}$, where we assume that the $f_i$ are locally Lipschitz for all $i=1,\ldots,n$. All signals ($x_i$ and inputs $u_i$, $i=1,\ldots,n$) are assumed to be right-continuous and to have left limits at all times and we denote $x_i^{-}(t):=\lim_{s\nearrow t}x_i(s)$, $u_i^{-}(t):=\lim_{s\nearrow t}u_i(s)$.

We define $N:=N_1+\ldots+N_n,\ M:=M_1+\ldots+M_n$, $x:=(x_1^T,\ldots,x_n^T)^T$, $u:=(u_1^T,\ldots,u_n^T)^T$, $f:=(f_1^T,\ldots,\ f_n^T)^T$ and $g:=(g_1^T,\ldots,\ g_n^T)^T$ such that the interconnected system \eqref{def:imp_sys_is} is of the form \eqref{def:imp_sys_single}. 
%the impulse time sequence of the whole system $\{t_k\}:=\left\{t|t=t_{i_{\tilde{k}}},\ \tilde{k}\in\N\right\},\ k\in\N$. 
%It may happen that at an impulse time $t_{i_{\tilde{k}}}$ there is a jump of the $i$th subsystem but not of the $j$th subsystem, $j\in\{1,...,n\},\ j\neq i$, for example. This circumstance may lead to a conservative condition for stability of the whole system, obtained from the exponential Lyapunov functions of the subsystems. Therefore, we define $I_k:=\{i|t_k=t_{i_{\tilde{k}}}\}$, which is the set of the impulse times of the $i$th subsystem and the whole system; $\overline{I}_k:=\{i|t_k\neq t_{i_{\tilde{k}}}\}\left|_{I_k}\right.$, which is the set of impulse times of the whole system, but not of the $i$th subsystem; and we denote $N(t,s)$ as the number of impulse times in the semi-open interval $(s,t]$ for the whole system. Then, at an impulse time $t_k$ of the whole system we define
%\begin{align*}
%	x\left|_{I_k}\right.&:=(0,\ldots,g_{i_1},\ldots,\ldots,g_{i_p},\ldots,0)^T, i_j\in I(t),\ j=1,\ldots,p,\\
%	x\left|_{\overline{I}_k}\right.&:=(0,\ldots,x_{i_1},\ldots,\ldots,x_{i_l},\ldots,0)^T,\ i_j\in \overline{I}(t),\ j=1,\ldots,l.
%\end{align*}
%
%With these definitions the interconnected system (\ref{def:imp_sys_is}) can be described as a system of the form \eqref{def:imp_sys_single}, where $g(x^-(t),u^{-}(t)):=x|_{I_k}+x|_{\overline{I}_k}$. 
We investigate under which conditions the whole system has the ISS property. In case of a system with several inputs the definition of ISS reads as follows:

Assume that a sequence $\{t_k\}$ is given. The $i$th subsystem of (\ref{def:imp_sys_single}) is \textit{ISS} if there exist $\beta_i\in\KL$, $\gamma_{ij}, \gamma_i \in{\cK}_\infty\cup\{$0$\}$ such that for every initial condition $x_i(t_0)$ and every input $u_i$ the corresponding solution to (\ref{def:imp_sys_is}) exists globally and satisfies  for all $t\geq t_0$
\begin{align}\label{eq:iss_is}
|x_i(t)|\leq \max\{\beta_i(|x_i(t_0)|,t-t_0),\max\limits_{j,j\neq i} \gamma_{ij}(\|x_j\|_{[t_0,t]}),\gamma_i(\|u\|_{[t_0,t]})\}. 
\end{align}
Functions $\gamma_{ij}$ are called gains. The impulsive system (\ref{def:imp_sys_is}) is {\it uniformly ISS} over a given class $\Simp$ of admissible sequences of impulse times if (\ref{eq:iss_is}) holds for every sequence in $\Simp$, with functions $\beta_i$ and $\gamma_i, \gamma_{ij}$ that are independent of the choice of the sequence.

Similarly, the Lyapunov functions for a system with several inputs are as follows: 
%\begin{definition}\label{def:iss_lyap_is}

%\makebox[-20pt][r]{\textbf{Kommentar!!!}}
%\fbox{Änderungen hier
%}

Assume that for each subsystem of the interconnected system \eqref{def:imp_sys_is} there is a given function $V_i:\R^{N_i}\rightarrow\R_+$, which is continuous, proper, positive definite and locally Lipschitz continuous on $\R^{N_i}\backslash\{0\}$.
For $i=1,\ldots,n$ the function $V_i$ is called an {\it exponential ISS-Lyapunov function} for the $i$th subsystem of \eqref{def:imp_sys_is} with {\it rate coefficients} $c_i, d_i\in \R$ if whenever $V_i(x_i)\geq\max\{\max\limits_{j, j\neq i}\gamma_{ij}(V_j(x_j)),\gamma_i(|u_i|)\}$ holds it follows
\begin{align}
\nabla V_i(x_i)\cdot f_i(x,u_i)\leq-c_iV_i(x_i)\ \text{f.a.a. }x,\text{ all }u_i\text{ and}\label{eq:iss_lyap1_is}
\end{align}
and for all $x$ and $u_i$ it holds 
\begin{align}
%V_i(g_i(x,u_i))\leq e^{-d_i}V_i(x_i)\ \text{for all }x, u_i, \label{eq:iss_lyap2_is}
V_i(g_i(x,u_i))\leq\max\{e^{-d_i}V_i(x_i),\max\limits_{j, j\neq i}\gamma_{ij}(V_j(x_j)),\gamma_i(|u_i|)\}, \label{eq:iss_lyap2_is}
\end{align}
where $\gamma_{ij}$, $\gamma_i$ are some functions from $\K_{\infty}$. A different formulation can be obtained by replacing \eqref{eq:iss_lyap2_is} by
\begin{align*}
	V_i(x_i)\geq\max\{\max\limits_{j, j\neq i}\tilde{\gamma}_{ij}(V_j(x_j)),\tilde{\gamma}_i(|u_i|)\}\ \Rightarrow\ V_i(g_i(x,u_i))\leq e^{-d_i}V_i(x_i),
\end{align*}
where $\tilde{\gamma}_{ij},\tilde{\gamma}_i\in\K_{\infty}$.
%\end{definition}

In general, even if all subsystems of (\ref{def:imp_sys_is}) are ISS, the whole system (\ref{def:imp_sys_single}) may be not ISS. Thus, we need conditions that guarantee ISS of (\ref{def:imp_sys_single}). In this paper we are concerned with an interconnection of impulsive systems that have exponential ISS-Lyapunov functions $V_i$ with linear gains $\gamma_{ij}$.

By slight abuse of notation we denote throughout the paper $\gamma_{ij}(r)= {\gamma}_{ij} r,\ \gamma_{ij},r \geq 0$.
To derive sufficient stability conditions for such an interconnection we collect the linear gains $\gamma_{ij}$ of the subsystems in a matrix $\Gamma=(\gamma_{ij})_{n\times n}$, $i,j=1,\dots,n$ denoting $\gamma_{ii} := 0$,\, $i=1,\dots,n$ for completeness, see \cite{DRW06b, DRW07, Rue07}. Note that this matrix describes in particular the interconnection topology of the whole network, moreover it contains the information about the mutual influence between the subsystems. We also introduce the gain operator $\Gamma:\R_{+}^n\rightarrow\R_{+}^n$ defined by
\begin{align}\label{operator_gamma}
\Gamma(s):=\left(\max_j\gamma_{1j} s_j,\ldots,\max_{j}\gamma_{nj}s_{j}\right)^T,\ s\in\R_{+}^n.
\end{align}

For the stability analysis of interconnected systems, we need the following:
\begin{definition}
We say that $\Gamma$ satisfies \textit{the small-gain condition}, if it satisfies
\begin{align}\label{smallgaincondition}
\Gamma(s)\not\geq s,\ \forall\ s\in\R^n_+\backslash\left\{0\right\}.
\end{align}
\end{definition}

As the gains in the matrix $\Gamma$ are linear, condition \eqref{smallgaincondition} is equivalent to 
\begin{align}\label{smallgaincondition_spectralradius}
\rho(\Gamma)<1,
\end{align}
where  $\rho$ is the spectral radius of $\Gamma$, see \cite{DRW06b, Rue07}. Note also that $\rho(\Gamma)<1$ implies that there exists a vector $s\in\R^n$, $s>0$ such that
\begin{align}\label{sigma cond 2 rho}
\Gamma(s)<s.
\end{align}

%
%To show one of the main results we need the notion of a so called $\Omega$-path, see \cite{DRW10,Rue10}.
%
%A function $\sigma=(\sigma_1,\dots,\sigma_n)^T:\R^n_+\rightarrow\R^n_+$, where $\sigma_i\in\K_\infty$ is called an {\it $\Omega$-path}, if it possesses the following properties:
%
%(i) $\sigma^{-1}_i$ is locally Lipschitz continuous on
%$(0,\infty)$;
%
%(ii) for every compact set $P\subset(0,\infty)$ there are finite
%constants $0<K_1<K_2$ such that for all points of differentiability of
%$\sigma^{-1}_i$ we have
%\begin{align*}%\label{sigma cond 1}
%0<K_1\leq(\sigma^{-1}_i)'(r)\leq K_2, \quad \forall r\in P
%\end{align*}
%(iii)%for all $r>0$
%\begin{align}\label{sigma cond 2}
%\Gamma(\sigma(r))<\sigma(r), \forall r>0.
%\end{align}
%The next theorem was proved in \cite{Rue10} and provides a condition for the existence of an $\Omega$-path.
%\begin{theorem}[Existence of $\Omega$-path]\label{Existenzsigma}
%Let $\Gamma\in({\cK}_{\infty}\cup\left\{0\right\})^{n\times n}$ be a gain matrix. If $\Gamma$ satisfies the small-gain condition
%\begin{align}\label{smallgaincondition}
%\Gamma(s)\not\geq s,\ \forall\ s\in\R^n_+\backslash\left\{0\right\},
%\end{align}
%then there exists an $\Omega$-path $\sigma$ with respect to $\Gamma$. This path can be chosen piecewise linear.
%\end{theorem}
%Note however that only the existence of $\Omega$-path is proved.

Now, we can formulate one of the main results that is an ISS small-gain theorem for impulsive systems without time-delays. This theorem allows to construct an exponential ISS-Lyapunov function for the whole interconnection.

\begin{theorem} \label{thm:main1}
%\textbf{Consider interconnection (\ref{def:imp_sys_is}) of systems (\ref{def:imp_sys_single}) with $t_{i_{\tilde{k}}}=t_k$.}
Assume that each subsystem of (\ref{def:imp_sys_is}) has an 
exponential ISS-Lyapunov function $V_i$ with corresponding linear ISS-Lyapunov
gains $\gamma_{ij}$ and rate coefficients $c_i$, $d_i$, $d_i\neq0$. 
If $\Gamma=(\gamma_{ij})_{n\times n}$ satisfies the small-gain condition \eqref{smallgaincondition}, then the exponential ISS-Lyapunov function for the whole system (\ref{def:imp_sys_is}) can be chosen as
\begin{align}\label{eq_Lyap}
	V(x):=\max_i\{\tfrac{1}{s_i}V_i(x_i)\},
\end{align}
where $s=(s_1,\ldots,s_n)^T$ is from \eqref{sigma cond 2 rho}. Its gain $\gamma$ is given by $\gamma(r):=\max\{e^{d},1\}\max_{i}\tfrac{1}{s_i}\gamma_i(r)$ and the rate coefficients are
$c:=\min\limits_{i}{c_i}$ and $d:=\min_{i,j,\ j\neq i}\{d_i,-\ln(\tfrac{s_j}{s_i}\gamma_{ij})\}$.
In particular, for all $\mu, \lambda>0$ (\ref{def:imp_sys_is}) is uniformly ISS over $\Simp[\mu,\lambda]$.
%
%Define $c:=\min\limits_{i}{c_i}$ and $d:=\min_{i,j,\ j\neq i}\{d_i,-\ln(\tfrac{s_j}{s_i}\gamma_{ij})\}$. For arbitrary constants $\mu, \lambda>0$, let $\Simp[\mu,\lambda]$ denote the class of impulse time sequences $\{t_k\}$ of the whole system. If the following holds
%
%i) $\Simp[\mu,\lambda]$ satisfies condition \eqref{eq:cond_int},
%
%ii) $\Gamma=(\gamma_{ij})_{n\times n}$ satisfies the small-gain condition \eqref{smallgaincondition}, 
%
%then the impulsive system (\ref{def:imp_sys_single}) is uniformly ISS over $\Simp[\mu,\lambda]$ and the exponential ISS-Lyapunov function is given by
%
% $\sigma=(\sigma_1,\ldots,\sigma_n)^T$ is a piecewise linear $\Omega$-path. 
\end{theorem}

%\begin{theorem}\label{thm:main1}
%%\textbf{Consider interconnection (\ref{def:imp_sys_is}) of systems (\ref{def:imp_sys_single}) with $t_{i_{\tilde{k}}}=t_k$.}
%Assume that each subsystem of (\ref{def:imp_sys_is}) has an 
%exponential ISS-Lyapunov function $V_i$ with corresponding linear ISS-Lyapunov
%gains $\gamma_{ij}$ and rate coefficients $c_i$, $d_i$, $d_i\neq0$. Define $c:=\min\limits_{i}{c_i}$ and $d:=\min_{i,j,\ j\neq i}\{d_i,-\ln(\tfrac{s_j}{s_i}\gamma_{ij})\}$. For arbitrary constants $\mu, \lambda>0$, let $\Simp[\mu,\lambda]$ denote the class of impulse time sequences $\{t_k\}$ of the whole system. If the following holds

%i) $\Simp[\mu,\lambda]$ satisfies condition \eqref{eq:cond_int},

%ii) $\Gamma=(\gamma_{ij})_{n\times n}$ satisfies the small-gain condition \eqref{smallgaincondition}, 

%then the impulsive system (\ref{def:imp_sys_single}) is uniformly ISS over $\Simp[\mu,\lambda]$ and the exponential ISS-Lyapunov function is given by
%\begin{align}\label{eq_Lyap}
	%V(x):=\max_i\{\tfrac{1}{s_i}V_i(x_i)\},
%\end{align}
%where $s=(s_1,\ldots,s_n)^T$ is from \eqref{sigma cond 2 rho}.
%% $\sigma=(\sigma_1,\ldots,\sigma_n)^T$ is a piecewise linear $\Omega$-path. 
%The gain $\gamma$ is given by $\gamma(r):=\max\{e^{d},1\}\max_{i}\tfrac{1}{s_i}\gamma_i(r).$
%\end{theorem}

The small-gain condition is used in \cite{DRW06b, DRW07, DRW10}, for example, to verify the ISS property of interconnected systems. Note that $\gamma_i$ are allowed to be nonlinear.

\textit{\textbf{Proof.}}
As the small gain condition \eqref{smallgaincondition} is satisfied, it follows from \eqref{smallgaincondition_spectralradius}
%Theorem~\ref{Existenzsigma} 
that there exists $s\in\R^n$, $s>0$ satisfying \eqref{sigma cond 2 rho}, see \cite{DRW06b}.
%an $\Omega$-path $\sigma$ with respect to $\Gamma$. We can choose this path to be piecewise linear. 
Let us  define $V$ as in \eqref{eq_Lyap}
%\begin{align}\notag
%V(x)=\max\limits_i\frac{1}{s_i}V_i(x_i)
%\end{align}
and show that this function is an exponential ISS-Lyapunov function for the system \eqref{def:imp_sys_single}. It can be easily checked that this function is locally Lipschitz, positive definite and radially unbounded.

For any $i\in\{1,\ldots,n\}$ consider open domains $M_i\in\R^N\setminus\{0\}$ defined by 
\begin{align}\label{eq:set_m}
\begin{aligned}
M_i:=&\{\left(x_1^T,\ldots,x_n^T\right)^T\in\R^N\backslash\left\{0\right\}:
%&\phantom{\{}
\tfrac{1}{s_i}V_i(x_i)>\max_{j\neq i}\tfrac{1}{s_j}V_j(x_j)\}.
\end{aligned}
\end{align}

Take arbitrary $\hat{x}=(\hat{x}_1^T,\ldots,\hat{x}_n^T)^T\in\R^N\backslash\{0\}$ for which there exist $i\in\{1,\ldots,n\}$, such that $\hat{x}\in M_i$. It follows that there is a neighborhood $U$ of $\hat{x}$ such that $V(x)=\frac{1}{s_i}V_i(x_i)$ is differentiable for almost all $x\in U$. 

We define $\bar{\gamma}(r):=\max_{i}\tfrac{1}{s_i}\gamma_i(r)$, $r>0$ and assume $V(x)\geq\bar{\gamma}(\left|u\right|)$. It follows from 
%\eqref{sigma cond 2}
\eqref{sigma cond 2 rho} that
\begin{align*}
	V_i(x_i)&=s_i V(x)\geq\max\{\max_{j} {\gamma}_{ij} s_j V(x),s_i\bar{\gamma}(|u|)\}\\
	&\geq\max\{\max_j {\gamma}_{ij} V_j(x_j),\gamma_i(|u_i|)\}.
\end{align*}
%Let $[s_l,s_{l+1}],\ s_l<s_{l+1},\ l=0,1,\ldots$ be an interval where $\sigma_i$ is linear, i.e., $\sigma_i(t)=a_{il}s,\ s\in[s_l,s_{l+1}]$. 
Then, 
%for all intervals with $\sigma_i(s)=a_{il}s$ and 
from \eqref{eq:iss_lyap1_is} we obtain for almost all $x$ %and $V_i(x_i)$%\in[s_l,s_{l+1}]$
\begin{align*}
	\dot{V}(x)=\tfrac{1}{s_i}\nabla V_i(x_i)\cdot f_i(x,u_i)	\leq-\tfrac{1}{s_i}c_iV_i(x_i)=-c_iV(x).
\end{align*}

We have shown that for $c=\min_ic_i$ the function $V$ satisfies \eqref{eq:iss_lyap1_ws} with $\bar{\gamma}$ for all $\hat{x} \in \cup_{i=1}^n M_i$. To treat the points $\hat{x} \in \R^N \backslash \cup_{i=1}^n M_i$ one can use the technique from \cite{DRW06b}.

%\textbf{Due to $t_{i_{\tilde{k}}}=t_k$ all the subsystems jump at the same time instants $t_k$. Thus $g(x^-(t),u^{-}(t))=x|_{I_k}$. } 
%The term $x|_{I_k}$ includes the jumps of the subsystems at an impulse time, whereas the term $x|_{\overline{I}_k}$ represents the subsystems with no jumps at an impulse time of the whole system, but with continuous behavior at time $t$, which was already considered in the calculations above to prove \eqref{eq:iss_lyap1_ws}. It remains to investigate the case $g(x^-(t),u^{-}(t))=x|_{I_k}$. Therefore, 
%\makebox[-20pt][r]{\textbf{Kommentar!!!}}
%\fbox{Änderungen hier
%}

With $d:=\min_{i,j,\ j\neq i}\{d_i,-\ln(\tfrac{s_j}{s_i}\gamma_{ij})\}$ and using \eqref{eq:iss_lyap2_is} it holds
\begin{align*}
V(g(x,u))=&\max_i\{\tfrac{1}{s_i}V_i(g_i(x_1,\ldots,x_n,u_i))\}\\
\leq&\max_i\{\tfrac{1}{s_i}\max\{e^{-d_i}V_i(x_i),\max\limits_{j, j\neq i}\gamma_{ij}V_j(x_j),\gamma_i(|u_i|)\}\}\\
\leq&\max_{i,j\ j\neq i}\{\tfrac{1}{s_i}e^{-d_i}s_iV(x),\tfrac{1}{s_i}\gamma_{ij}s_j V(x),\tfrac{1}{s_i}\gamma_i(|u_i|)\}\\
\leq&\max\{e^{-d}V(x),\bar{\gamma}(|u|)\}.
\end{align*}
%$V_i(x_i)$\in[s_l,s_{l+1}]$,

Define $\tilde{\gamma}(r):=e^{d}\bar{\gamma}(r)$. If it holds $V(x)\geq\tilde{\gamma}(|u|)$, it follows
\begin{align*}
	V(g(x,u))\leq\max\{e^{-d}V(x),\bar{\gamma}(|u|)\}=\max\{e^{-d}V(x),e^{-d}\tilde{\gamma}(|u|)\}\leq e^{-d}V(x)
\end{align*}
for all $x$ and $u$ and $V$ satisfies \eqref{eq:iss_lyap2_ws} with $\tilde{\gamma}$. Define $\gamma(r):=\max\{e^{d},1\}\max_{i}\tfrac{1}{s_i}\gamma_i(r)$, then we conclude that $V$ is an exponential ISS-Lyapunov function with rate coefficients $c, d$ and gain $\gamma$.
%satisfies condition (\ref{eq:iss_lyap2_ws}). 

All conditions of Definition~\ref{def:iss_lyap_ws} are satisfied and thus $V$ is the exponential ISS-Lyapunov function of the system \eqref{def:imp_sys_single}. We can apply Theorem~\ref{thm:stab_1sys} and the overall system is uniformly ISS over $S[\mu,\lambda]$, $\mu,\lambda>0$.
\hfill $\Box$

%\makebox[-20pt][r]{\textbf{Kommentar!!!}}
%\fbox{Geändert hier.
%}

\begin{remark}
Using $c=\min\limits_{i} c_i,\ d=\min_{i,j,\ j\neq i}\{d_i,-\ln(\tfrac{s_j}{s_i}\gamma_{ij})\}$ in the dwell-time condition some kind of conservativeness may occur, which means that the ISS property for an interconnected impulsive system cannot be verified by the application of Theorem~\ref{thm:main1}, although the system possesses the ISS property.
\end{remark}

\subsection{Example with nonlinear gains}

We have formulated Theorem~\ref{thm:main1} for the case of linear gains. Note that not for all interconnections with nonlinear gains one can construct an exponential Lyapunov function for the whole system, even if the small-gain condition is satisfied. However, if the small-gain condition holds, then we can construct non-exponential Lyapunov functions for the whole system as in \cite{DRW06b}, \cite{DRW10}. In this case the dwell-time condition from \cite{HLT08} cannot be applied and one has to develop more general conditions that guarantee ISS of the system.

Nevertheless, in some special cases the treatment of the interconnections with nonlinear gains is also possible.
For example, if the continuous and discrete dynamics of all subsystems are stabilizing,
i.e., with $c_i, d_i>0$,
and the small-gain condition holds,
 then we can construct an ISS-Lyapunov function (non-exponential in general) for the whole system using the methodology from \cite{DRW10}. Then, according to Theorem~2 in \cite{HLT08} the interconnection will be ISS.
Note that in this case the dwell-time condition is not needed anymore. 

For the case that the discrete or continuous dynamics is unstable we consider the following example with a construction of an ISS-Lyapunov function.

Let $T = \{t_k\}$ be a sequence of impulse times. Consider two interconnected nonlinear systems
\begin{align*}
\dot{x}_1(t) =&\ -x_1(t) + x_2^2(t), \ t \notin T,\\
x_1(t)=&\ e x_1^-(t),  \ t \in T
\end{align*}
and
\begin{align*}
\dot{x}_2(t) =&\ -x_2(t) + \frac{1}{2} \sqrt{|x_1(t)|},  \ t \notin T,\\
x_2(t)=&\ e x_2^-(t),  \ t \in T.
\end{align*}
The exponential ISS-Lyapunov functions and Lyapunov gains for these subsystems are given by 
\begin{align*}
V_1(x_1)=|x_1|, \quad   \gamma_{12}(r)=\tfrac{1}{1-\eps_1} r^2,\\
V_2(x_2)=|x_2|, \quad   \gamma_{21}(r)=\tfrac{1}{2(1-\eps_2)} \sqrt{r},
\end{align*}
where $\eps_1$, $\eps_2 \in (0,1)$. We have the estimates
\begin{align*}
|x_1| \geq \gamma_{12} (|x_2|) \Rightarrow \dot{V}_1(x_1) \leq - \eps_1 V_1(x_1),\\
|x_2| \geq \gamma_{21} (|x_1|) \Rightarrow \dot{V}_2(x_2) \leq - \eps_2 V_2(x_2).
\end{align*}
%Thus, in the dwell- time condition \eqref{Dwell-Time-Cond}
%$c_1=c_2=\eps$ and $d_1 = d_2 = -1$, and $\eps$ can be made close to one.
The small-gain condition 
\begin{align}
\label{KGB_NichtLin}
\gamma_{12} \circ \gamma_{21} (r) = \frac{1}{4(1-\eps_1)(1-\eps_2)^2} r < r, \quad \forall r>0
\end{align}
is satisfied, if 
\begin{align}
\label{KGB_UnserSyst}
h(\eps_1,\eps_2):=(1-\eps_1)(1-\eps_2)^2 > \frac{1}{4}
\end{align}
holds.
%The Lyapunov function for interconnection can be constructed by using of so-called $\Omega$-path $\sigma=(\sigma_1,\sigma_2)$ (see \cite{DRW06b}). In our case it can be chosen as 
%\[
%\sigma_1(r)=r,\quad  \sigma_2(r)=a\sqrt{r},\ \forall r, \quad a \in (\frac{1}{2(1-\eps_2)},\sqrt{1-\eps_1}).
 %\]
%Then
%\[
%\sigma_2^{-1}(r)=\frac{1}{a^2} r^2,\ \forall r \geq 0.
%\]
In this case, an ISS-Lyapunov function for the interconnection is given by
\begin{align*}
V(x)=\max\{ |x_1|, \tfrac{1}{a^2} x_2^2  \} , \quad \text{where }a \in (\tfrac{1}{2(1-\eps_2)},\sqrt{1-\eps_1})
\end{align*}
and we have the estimate
 \begin{align}
\label{LyapAbsch_Diskr}
V(g(x))=V(e\cdot x) \leq e^2 V(x).
\end{align}
Thus, $d=-2$ for the interconnection. The Lyapunov estimates of the continuous dynamics for $V$ are as follows:
\begin{align*}
&\frac{d}{dt} (|x_1|) \leq -\eps_1  |x_1|, \\
&\frac{d}{dt} \left(\frac{1}{a^2} x_2^2 \right) = \frac{d}{dt} \left(\frac{1}{a^2} V_2(x_2)^2\right) \leq 
-2 \eps_2 \frac{1}{a^2} \left( V_2(x_2) \right)^2. 
\end{align*}
Using estimates as in the proof of Theorem~\ref{thm:main1}, we obtain
%\ref{KleinGainSatz_ImpSys}:
 \begin{align}
\label{LyapAbsch_Cont}
\frac{d}{dt} V(x) \leq - \min\{\eps_1, 2\eps_2 \} V(x).
\end{align}

Function $h$, defined by \eqref{KGB_UnserSyst}, is increasing on both arguments (remember, that $\eps_i \in (0,1)$), which implies that to maximize $\eps:= \min\{\eps_1, 2\eps_2 \}$, we have to choose $\eps_1 = 2 \eps_2$. Then, from \eqref{KGB_NichtLin} we obtain the inequality
\begin{align*}
(1-2\eps_2)(1-\eps_2)^2 > \frac{1}{4}.
\end{align*}
Thus, the best choice for $\eps_2$ is $\eps_2 \approx 0.267$. From the inequalities \eqref{LyapAbsch_Diskr} and \eqref{LyapAbsch_Cont} we see that if the dwell-time condition \eqref{eq:cond_int}
%\eqref{Dwell-Time-Cond} 
with $d=-2$ and $c=2 \cdot 0.267$ is satisfied, then the system under investigation is stable according to Theorem 1 in \cite{HLT08}.

In this example we were able to construct an exponential ISS Lyapunov function for the whole system due the special type of nonlinearities in internal gains. For general nonlinear internal gains this construction is not always possible.

%%%%%%%%%%%%%%%%%%%%%%%%%%%%%%%%%%%%%%%%%%%%%%%%%%%%%%%%%%%%%%%%%%%%%%%%%%%%%%%%%%%%%%%%%%%%%%%%%%%%%%%%%%%%%%%%%%%%%%%%%%%%%%%%%%%%%%%%%%%%%%%%%%%%%%%%%%%%%%%%%%%
\subsection{Systems with time-delays}
\label{SubSec:ZeitVerz}

Now we consider $n$ interconnected impulsive systems with time-delays of the form
\begin{align}\label{Gekdelaysys}
\begin{aligned}
	\dot{x}_i(t)&=f_i(x_1^t,\ldots,x_n^t,u_i(t)),\ t\neq t_{k},\\
	x_i(t)&=g_i((x_1^t)^-,\ldots,(x_n^t)^-,u_i^-(t)),\ t=t_{k},
	\end{aligned}
\end{align}
where the same assumptions on the system as in the delay-free case are considered with the following differences: We denote $x^t_i(\tau):=x_i(t+\tau),\ \tau\in\left[-\theta,0\right]$, where $\theta$ is the maximum involved delay and $(x_i^t)^{-}(\tau):=x^-_i(t+\tau):=\lim_{s\nearrow t}x_i(s+\tau)$, $\tau\in[-\theta,0]$. Furthermore, $f_i:PC([-\theta,0],\R^{N_1})\times\ldots\times PC([-\theta,0],\R^{N_n})\times\R^{M_i}\rightarrow\R^{N_i}$, and $g_i:PC([-\theta,0],\R^{N_1})\times\ldots\times PC([-\theta,0],\R^{N_n})\times\R^{M_i}\rightarrow\R^{N_i}$, where we assume that $f_i$ are locally Lipschitz, $i=1,\ldots,n$.

If we define $\ N,\ M,\ x,\ u,\ f$ and $g$ as in the delay-free case, then (\ref{Gekdelaysys}) becomes the system of the form \eqref{eq:SingleIMPTDS}. The ISS property for systems with several inputs and time-delays is as follows:

Suppose that a sequence $\{t_k\}$ is given. The $i$th subsystem of (\ref{Gekdelaysys}) is \textit{ISS}, if there exist $\beta_i\in\KL$, $\gamma_{ij}, \gamma_i^u \in {\cK}_\infty\cup\{0\}$ such that for every initial condition $\xi_i$ and every input $u_i$ the corresponding solution to the $i$th subsystem of (\ref{Gekdelaysys}) exists globally and satisfies
\begin{align}\label{eq:iss_is_delay}
|x_i(t)|\leq \max\{\beta_i(\|\xi_i\|_{[-\theta,0]},t-t_0),\max\limits_{j,j\neq i} \gamma_{ij}(\|x_j\|_{[t_0-\theta,t]}),\gamma^u_i(\|u\|_{[t_0,t]})\} 
\end{align}
for all $t\geq t_0$. The $i$th subsystem of (\ref{Gekdelaysys}) is {\it uniformly ISS} over a given class $\Simp$ of admissible sequences of impulse times if (\ref{eq:iss_is_delay}) holds for every sequence in $\Simp$, with functions $\beta_i,\ \gamma_{ij}$ and $\gamma^u_i$ that are independent of the choice of the sequence.

In the following subsections, we present useful tools to analyze a system of the form (\ref{Gekdelaysys}) in view of stability: ISS-Lyapunov-Razumikhin functions and ISS-Lyapunov-Krasovskii functionals for the subsystems.
%%%%%%%%%%%%%%%%%%%%%%%%%%%%%%%%%%%%%%%%%%%%%%%%%%%%%%%%%%%%%%%%%%%%%%%%%%%%%%%%%%%%%%%%%%%%%%%%%%%%%%%%%%%%%%%%%%%%%%%%%%%%%%%%%%%%%%%%%%%%%%%%%%%%%%%%%%%%%%%%%%%
\subsubsection{Lyapunov-Razumikhin functions}
%The definition of ISS-Lyapunov-Razumikhin functions for impulsive systems with several inputs and time-delays reads as follows:
%\begin{definition}\label{DefLRgek}

%\makebox[-20pt][r]{\textbf{Kommentar!!!}}
%\fbox{Änderungen hier
%}

Assume that for each subsystem of the interconnected system \eqref{Gekdelaysys} there is a given function $V_i:\R^{N_i}\rightarrow\R_+$, which is continuous, proper, positive definite and locally Lipschitz continuous on $\R^{N_i}\backslash\{0\}$.
For $i=1,\ldots,n$ the function $V_i$ is called an \textsl{exponential ISS-Lyapunov-Razumikhin function of the $i$th subsystem of (\ref{Gekdelaysys})}, if there exist $\gamma_i^u\in{\cK}\cup\{0\},\ \gamma_{ij}\in{\cK}_{\infty}\cup\{0\},\ j=1,\ldots,n$ and scalars $c_i,d_i\in\R,$ such that whenever $V_i(\phi_i(0))\geq \max\{\max_{j}\gamma_{ij}(\|V_j^t(\phi_j)\|),\gamma_i^u(|u_i|)\}$ holds it follows
\begin{align}
	\D^+V_i(\phi_i(0))\leq-c_iV_i(\phi_i(0))\label{CondLR2gek}
	\end{align}
for all $\phi=(\phi_1^T,\ldots,\phi_n^T)^T\in PC([-\theta,0],\R^{N})$ and $u_i\in\R^{M_i}$, where $V_j^t:PC\left(\left[-\theta,0\right];\R^{N_j}\right) \to PC\left(\left[-\theta,0\right];\R^{N_j}\right)$ with $V_j^t(x_j^t)(\tau):=V_j(x_j(t+\tau)),\ \tau\in\left[-\theta,0\right]$ and it holds
	\begin{align}	V_i(g_i(\phi_1,\ldots,\phi_n,u_i))\leq\max\{e^{-d_i}V_i(\phi_i(0)),\max_{j}\gamma_{ij}(\|V_j^t(\phi_j)\|),\gamma_i^u(|u_i|)\},\label{CondLR3gek}
\end{align}
for all $\phi\in PC([-\theta,0],\R^{N})$ and $u_i\in\R^{M_i}$. Another formulation can be obtained by replacing \eqref{CondLR3gek} by
\begin{align*}
	V_i(\phi_i(0))\geq \max\{\max_{j}\tilde{\gamma}_{ij}(\|V_j^t(\phi_j)\|),\tilde{\gamma}_i^u(|u_i|)\}\ \Rightarrow\ V_i(g_i(\phi,u_i))\leq e^{-d_i}V_i(\phi_i(0)),
\end{align*}
where $\tilde{\gamma}_{ij},\tilde{\gamma}_i^u\in\K_{\infty}$.

We consider linear gains $\gamma_{ij}$ and define the \textsl{gain-matrix} $\Gamma:=(\gamma_{ij})_{n\times n}$ and the map $\Gamma:\R^n_+\rightarrow\R^n_+$ by $\Gamma(s):=\left(\max_{j}\gamma_{1j}s_j,\ldots,\max_{j}\gamma_{nj}s_j\right)^T,\ s\in\R^n_+.$ 
%\end{definition}

Now, we state one of our main results: the ISS small-gain theorem for interconnected impulsive systems with time-delays and linear gains $\gamma_{ij}$. We construct the Lyapunov-Razumikhin function and the gain of the overall system under a small-gain condition, which is here of the form 
\begin{align}\label{eq:sgcequivalentLR}
\begin{aligned}
	&\ \Gamma(s)\not\geq \min\{e^{-\mu},e^{-d-\mu}\}s, \forall\ s\in\R^n_+\backslash\left\{0\right\}\\
	\Leftrightarrow&\ \exists s\in\R^n_+\backslash\left\{0\right\}:\ \Gamma(s)<\min\{e^{-\mu},e^{-d-\mu}\}s,
	\end{aligned}
\end{align}
where $\mu$ is from the dwell-time condition \eqref{eq:cond_int} and $d:=\min_i d_i$. The dwell-time condition on the size of the time intervals between impulses is the same as in the delay-free case.

\begin{theorem}\label{th:main2}
%\textbf{Consider interconnection (\ref{eq:SingleIMPTDS}) of systems (\ref{Gekdelaysys}) with $t_{i_{\tilde{k}}}=t_k$.}
Assume that each subsystem of (\ref{Gekdelaysys}) has an exponential ISS-Lyapunov-Razumikhin function with $c_i,d_i\in\R,\ d_i\neq0$ and gains $\gamma_i^u$, $\gamma_{ij}$, where $\gamma_{ij}$ are linear.
Define $c:=\min_ic_i$ and $d:=\min_{i}d_i$. %For arbitrary constants $\mu, \lambda>0$, let $\Simp[\mu,\lambda]$ denote the class of impulse time sequences $\{t_k\}$ of the whole system. 
If for some $\mu>0$ the operator $\Gamma$ satisfies the small-gain condition \eqref{eq:sgcequivalentLR},
%\begin{align}\label{sigma cond 2 rhoNEU}
	%\Gamma(s)\not\geq \min\{e^{-\mu},e^{-d-\mu}\}s, \forall\ s\in\R^n_+\backslash\left\{0\right\}
%\end{align}
then for all $\lambda>0$  the whole system \eqref{eq:SingleIMPTDS} is uniformly ISS over $\Simp[\mu,\lambda]$ and the exponential ISS-Lyapunov-Razumikhin function is given by $V(x):=\max_i\{\tfrac{1}{s_i}V_i(x_i)\}$, where $s=(s_1,\ldots,s_n)^T$ is from \eqref{eq:sgcequivalentLR}.
% $\sigma=(\sigma_1,\ldots,\sigma_n)^T$ is a piecewise linear $\Omega$-path. 
The gains are given by $\gamma_t(r):=\max\{e^{d},1\}\max_{k,j}\frac{1}{s_k}{\gamma}_{kj}s_j r,$ $\gamma_u(r):=\max\{e^{d},1\}\max_{i}\tfrac{1}{s_i}\gamma_i^u(r).$
\end{theorem}

%\begin{theorem}\label{th:main2}
%%\textbf{Consider interconnection (\ref{eq:SingleIMPTDS}) of systems (\ref{Gekdelaysys}) with $t_{i_{\tilde{k}}}=t_k$.}
%Assume that each subsystem of (\ref{Gekdelaysys}) has an exponential ISS-Lyapunov-Razumikhin function with $c_i,d_i\in\R,\ d_i\neq0$ and gains $\gamma_i^u$, $\gamma_{ij}$, where $\gamma_{ij}$ are linear.
%Define $c:=\min_ic_i$ and $d:=\min_{i}d_i$. %For arbitrary constants $\mu, \lambda>0$, let $\Simp[\mu,\lambda]$ denote the class of impulse time sequences $\{t_k\}$ of the whole system. 
%If $\Gamma$ satisfies the small-gain condition 
%\begin{align}\label{sigma cond 2 rhoNEU}
	%\Gamma(s)\not\geq \min\{e^{-\mu},e^{-d-\mu}\}s, \forall\ s\in\R^n_+\backslash\left\{0\right\}
%\end{align}
%then the whole system \eqref{eq:SingleIMPTDS} is uniformly ISS over $\Simp[\mu,\lambda]$, $\mu,\lambda>0$ and the exponential ISS-Lyapunov-Razumikhin function is given by $V(x):=\max_i\{\tfrac{1}{s_i}V_i(x_i)\}$, where $s=(s_1,\ldots,s_n)^T$ is from \eqref{eq:sgcequivalentLR}.
%% $\sigma=(\sigma_1,\ldots,\sigma_n)^T$ is a piecewise linear $\Omega$-path. 
%The gains are given by $\gamma_t(r):=\max\{e^{d},1\}\max_{k,j}\frac{1}{s_k}{\gamma}_{kj}s_j r,$ $\gamma_u(r):=\max\{e^{d},1\}\max_{i}\tfrac{1}{s_i}\gamma_i^u(r).$
%\end{theorem}

The proof goes along the lines of the proof of Theorem~\ref{thm:main1} with corresponding changes due to time-delay systems and the additional gain $\gamma_t(r)$. 

\textit{\textbf{Proof.}}

%Let $0\neq x=(x_1^T,\ldots,x_n^T)^T$. 
We define $V(x)$ as in \eqref{eq_Lyap} and show that $V$ is the exponential ISS-Lyapunov-Razumikhin function for the overall system. Note that $V$ is locally Lipschitz continuous, positive definite and radially unbounded. For any $i\in\{1,\ldots,n\}$ consider open domains $M_i\in\R^N\backslash\{0\}$ defined as in \eqref{eq:set_m}.

Take arbitrary $\hat{x}=(\hat{x}_1^T,\ldots,\hat{x}_n^T)^T\in\R^N\backslash\{0\}$ for which there exist $i\in\{1,\ldots,n\}$, such that $\hat{x}\in M_i$. It follows that there is a neighborhood $U$ of $\hat{x}$ such that $V(x)=\frac{1}{s_i}V_i(x_i)$ and $D^+ V(x)$ exists. 

%Now for any $\hat{x}=(\hat{x}_1^T,\ldots,\hat{x}_n^T)^T\in\R^N\backslash\{0\}$ there is at least one $i\in\{1,\ldots,n\}$ such that $\hat{x}\in M_i$ and it follows that there is a neighborhood $U$ of $\hat{x}$ such that $V(x)=\tfrac{1}{s_i}V_i(x_i)$ and $D^+ V(x)$ exists. 

We define the gains $\bar{\gamma}_t(r):=\max_{k,j}\frac{1}{s_k}{\gamma}_{kj} s_jr$, $\bar{\gamma}_u(r):=\max_{i}\tfrac{1}{s_i}\gamma_i^u(r)$, $r>0$ and assume $V(x(t))\geq\max\{\bar{\gamma}_t(\|V^t(x^t)\|),\bar{\gamma}_u(\left|u(t)\right|)\}$. %Note that $\bar{\gamma}_d(r)<\min\{e^{-\mu},e^{d-\mu}\}r$, by \eqref{eq:sgcequivalentLR}. 
It follows
\begin{align*}
	V_i(x_i(t))\geq&\ s_i\max\{\max_{k,j} \frac{1}{s_k} {\gamma}_{kj} s_j \|V^t(x^t)\|,
%	\\	&&\phantom{s_i\cdot\max\{}
	\max_{i}\tfrac{1}{s_i}\gamma_i^u(|u(t)|)\}\\
	\geq&\ \max\{\max_j {\gamma}_{ij}\|V_j^t(x_j^t)\|,\gamma_i^u(|u_i(t)|)\}.
\end{align*}
%Let $[s_l,s_{l+1}],\ s_l<s_{l+1},\ l=0,1,\ldots$ be an interval where $\sigma_i$ is linear, i.e., $\sigma_i(s)=a_{il}s,\ s\in[s_l,s_{l+1}]$. 
Then, %for all intervals with $\sigma_i(s)=a_{il}s$ and 
from (\ref{CondLR2gek}) we obtain
\begin{align*}
	\D^+V(x(t))=\D^+\tfrac{1}{s_i}V_i(x_i(t))\leq-\tfrac{1}{s_i}c_iV_i(x_i(t))=-c_iV(x(t)).
\end{align*}

We have shown that for $c=\min_ic_i$ the function $V$ satisfies (\ref{CondLR2}) with $\bar{\gamma}_t,\ \bar{\gamma}_u$ 
for all $\hat{x} \in \cup_{i=1}^n M_i$. To treat the points $\hat{x} \in \R^N \backslash \cup_{i=1}^n M_i$ one can use the technique from \cite{DRW06b} or \cite{DaM11}.

%By definition of $c=\min_ic_i$ the function $V$ satisfies (\ref{CondLR2}) with $\bar{\gamma}_t,\ \bar{\gamma}_u$. 

With $d:=\min_{i}d_i$ and using \eqref{CondLR3gek} it holds
\begin{align*}
V(g(x^t,u(t)))=&\max_i\{\tfrac{1}{s_i}V_i(g_i(x_1^t,\ldots,x_n^t,u_i(t)))\}\\
\leq&\max_i\{\tfrac{1}{s_i}\max\{e^{-d_i}V_i(x_i(t)),\max_{j}\gamma_{ij}\|V_j^t(x_j^t)\|,\gamma_i^u(|u_i(t)|)\}\}\\
\leq&\max_{i,j\ j\neq i}\{\tfrac{1}{s_i}e^{-d_i}s_iV(x(t)),\tfrac{1}{s_i}\gamma_{ij}s_j \|V^t(x^t)\|,\tfrac{1}{s_i}\gamma_i^u(|u_i(t)|)\}\\
\leq&\max\{e^{-d}V(x(t)),\bar{\gamma}_t(\|V^t(x^t)\|),\bar{\gamma}_u(|u(t)|)\}.
\end{align*}
%$V_i(x_i)$\in[s_l,s_{l+1}]$,
Define $\tilde{\gamma}_t(r):=e^{d}\bar{\gamma}_t(r)$ and $\tilde{\gamma}_u(r):=e^{d}\bar{\gamma}_u(r)$.

If $V(x(t))\geq\max\{\tilde{\gamma}_t(\|V^t(x^t)\|),\tilde{\gamma}_u(|u(t)|)\}$ holds, it follows
\begin{align*}	
V(g(x^t,u(t)))&\leq\max\{e^{-d}V(x(t)),\bar{\gamma}_t(\|V^t(x^t)\|),\bar{\gamma}_u(|u(t)|)\}\\
&=\max\{e^{-d}V(x(t)),e^{-d}\tilde{\gamma}_t(\|V^t(x^t)\|),e^{-d}\tilde{\gamma}_u(|u(t)|)\}\\
	&\leq e^{-d}V(x(t)),
\end{align*}
i.e., $V$ satisfies the condition \eqref{CondLR3} with $\tilde{\gamma}_t,\ \tilde{\gamma}_u$. Now, define $\gamma_t(r):=\max\{\bar{\gamma}_t(r),\tilde{\gamma}_t(r)\}$ and $\gamma_u(r):=\max\{\bar{\gamma}_u(r),\tilde{\gamma}_u(r)\}$. Then, $V$ satisfies (\ref{CondLR2}) and \eqref{CondLR3} with $\gamma_t,\ \gamma_u$.

By \eqref{eq:sgcequivalentLR} it holds \\
$\gamma_t(r)=\max\{\bar{\gamma}_t(r),e^{d}\bar{\gamma}_t(r)\}<
\max\{\min\{e^{-\mu},e^{-d-\mu}\}, \min\{e^{d-\mu},e^{-\mu}\}\}r = e^{-\mu}r$.
%$\gamma_t(r)=\max\{\bar{\gamma}_t(r),e^{d}\bar{\gamma}_t(r)\}<\min\{e^{-\mu},e^{d-\mu},e^{-d-\mu}\}r<e^{-\mu}r$.

All conditions of Definition~\ref{DefLRfunction} are satisfied and $V$ is the exponential ISS-Lyapunov-Razumikhin function of the whole system of the form \eqref{eq:SingleIMPTDS}. We can apply Theorem~\ref{PropISSRazumikhin} and the whole system is uniformly ISS over $\Simp [\mu,\lambda]$.\hfill $\Box$

%\begin{remark}
%The meaning of $c,d\in\R$ is the same as in the delay-free case. If both dynamics, the continuous and the discontinuous, are stable, which means that $c,d>0$, then the condition \eqref{eq:cond_int} in Theorem~\ref{th:main2} is not necessary, according to Theorem~2 in \cite{HLT08}. 
%\end{remark}In the next subsection we prove a similar result as in Theorem~\ref{th:main2} but with Lyapunov-Krasovskii functionals.

%%%%%%%%%%%%%%%%%%%%%%%%%%%%%%%%%%%%%%%%%%%%%%%%%%%%%%%%%%%%%%%%%%%%%%%%%%%%%%%%%%%%%%%%%%%%%%%%%%%%%%%%%%%%%%%%%%%%%%%%%%%%%%%%%%%%%%%%%%%%%%%%%%%%%%%%%%%%%%%%%%%
\subsubsection{Lyapunov-Krasovskii functionals}
Let us consider $n$ interconnected impulsive subsystems of the form
\begin{align}\label{eq:InterIMPTDS}
\begin{aligned}
	\dot{x}_i(t)=&f_i(x_1^t,\ldots,x_n^t,u_i(t)),\ t\neq t_k,\\
	x_i^t=&g_i((x_1^t)^-,\ldots,(x_n^t)^-,u_i^-(t)),\ t=t_k,
	\end{aligned}
\end{align}
$k\in\N,\ i=1,\ldots,n$, where we make the same assumptions as in the previous subsections and the functionals $g_i$ are now maps from $PC\left(\left[-\theta,0\right];\R^{N_1}\right)\times\ldots\times PC\left(\left[-\theta,0\right];\R^{N_n}\right)\times\R^{M_i}$ into $PC\left(\left[-\theta,0\right];\R^{N_i}\right)$. 

%As in the Subsection~\ref{subsec:LKfunctionalsingle} the current state and the current input together determine the current output and the next state of the system. Therefore, we change the discontinuous behavior in \eqref{eq:InterIMPTDS} and this causes that the Lyapunov-Razumikhin approach cannot be applied. We use Lyapunov-Krasovskii functionals for the stability analysis of systems of the form \eqref{eq:InterIMPTDS}.
Note that the ISS property for systems of the form \eqref{eq:InterIMPTDS} is the same as in \eqref{eq:iss_is_delay}.

If we define $N,\ M,\ x,\ u,\ f$ and $g$ as in the previous subsection, then (\ref{eq:InterIMPTDS}) becomes the system of the form \eqref{eq:SingleIMPTDSstatejump}. ISS-Lyapunov-Krasovskii functionals of systems with several inputs and time-delays are as follows:
%\begin{definition}\label{DefLKgek}

%\makebox[-20pt][r]{\textbf{Kommentar!!!}}
%\fbox{Änderungen hier
%}

Assume that for each subsystem of the interconnected system \eqref{eq:InterIMPTDS} there is a given functional $V_i:PC\left(\left[-\theta,0\right];\R^{N_i}\right)\rightarrow\R_+$, which is locally Lipschitz continuous, positive definite and radially unbounded.
For $i=1,\ldots,n$ the functional $V_i$ is called an \textsl{exponential ISS-Lyapunov-Krasovskii functional of the $i$th subsystem of \eqref{eq:InterIMPTDS}}, if there exist $\gamma_i\in{\cK}\cup\{0\},\ \gamma_{ij}\in{\cK}_{\infty}\cup\{0\},\ \gamma_{ii}\equiv0,\ i,j=1,\ldots,n$ and scalars $c_i,d_i\in\R$ such that whenever $V_i(\phi_i)\geq \max\{\max_{j}\gamma_{ij}(V_j(\phi_j)),\gamma_i(|u_i|)\}$ holds it follows
\begin{align}
\D^+V_i(\phi_i,u_i)\leq-c_iV_i(\phi_i)\label{CondLK2gek}
\end{align}
and 
\begin{align}
%V_i(g_i(x,u_i))\leq e^{-d_i}V_i(x_i)\ \text{for all }x, u_i, \label{eq:iss_lyap2_is}
V_i(g_i(\phi,u_i))\leq\max\{e^{-d_i}V_i(\phi_i),\max\limits_{j, j\neq i}\gamma_{ij}(V_j(\phi_j)),\gamma_i(|u_i|)\}, \label{CondLK3gek}
\end{align}
for all $\phi\in PC\left(\left[-\theta,0\right];\R^{N_i}\right)$, $u_i\in\R^{M_i}$. A different formulation can be obtained by replacing \eqref{CondLK3gek} by
\begin{align*}
	V_i(\phi_i)\geq\max\{\max\limits_{j, j\neq i}\tilde{\gamma}_{ij}(V_j(\phi_j)),\tilde{\gamma}_i(|u_i|)\}\ \Rightarrow\ V_i(g_i(\phi,u_i))\leq e^{-d_i}V_i(\phi_i),
\end{align*}
where $\tilde{\gamma}_{ij},\tilde{\gamma}_i\in\K_{\infty}$.

%V_i(g_i(\phi,u_i))\leq e^{-d_i}V_i(\phi_i),
%\end{align}
%for all $\phi_i\in PC\left(\left[-\theta,0\right];\R^{N_i}\right),\ \phi=(\phi_1^T,\ldots,\phi_n^T)^T$ and $u_i\in\R^{M_i}$. Furthermore, we define the \textsl{gain-matrix} $\Gamma:=(\gamma_{ij})_{n\times n}$ with $\gamma_{ii}\equiv0$ and the map $\Gamma:\R^n_+\rightarrow\R^n_+$ by $\Gamma(s):=(\max_{j}\gamma_{1j}(s_j),\ldots,\max_{j}\gamma_{nj}(s_j))^T,\ s\in\R^n_+$. 
%\end{definition}

%\makebox[-20pt][r]{\textbf{Kommentar!!!}}
%\fbox{Änderungen hier nur bei der Definition von $d$
%}

The next result is an ISS small-gain theorem for impulsive systems with time-delays using the Lyapunov-Krasovskii methodology. This theorem allows to construct an exponential ISS-Lyapunov-Krasovskii functional and the corresponding gain for the whole interconnection under a dwell-time and a small-gain condition.
\begin{theorem}\label{thm:main3}
%\textbf{Consider interconnection (\ref{eq:SingleIMPTDSstatejump}) of systems (\ref{eq:InterIMPTDS}) with $t_{i_{\tilde{k}}}=t_k$.}
Assume that each subsystem of \eqref{eq:InterIMPTDS} has an exponential ISS-Lyapunov-Krasovskii functional $V_i$ with corresponding gains $\gamma_i,\gamma_{ij}$, where $\gamma_{ij}$ are linear, and rate coefficients $c_i$, $d_i$, $d_i\neq0$. Define $c:=\min\limits_{i}{c_i}$ and $d:=\min_{i,j,\ j\neq i}\{d_i,-\ln(\tfrac{s_j}{s_i}\gamma_{ij})\}$. %For arbitrary constants $\mu, \lambda>0$, let $\Simp[\mu,\lambda]$ denote the class of impulse time sequences $\{t_k\}$. 
If $\Gamma=(\gamma_{ij})_{n\times n}$, $\gamma_{ii}\equiv0$ satisfies the small-gain condition \eqref{smallgaincondition}, then the impulsive system \eqref{eq:SingleIMPTDSstatejump} is uniformly ISS over $\Simp[\mu,\lambda]$, $\mu,\lambda>0$ and the exponential ISS-Lyapunov-Krasovskii functional is given by
\begin{align}
	V(\phi):=\max_i\{\tfrac{1}{s_i}V_i(\phi_i)\},
\end{align}
where $s=(s_1,\ldots,s_n)^T$ is from \eqref{sigma cond 2 rho}, $\phi\in PC\left(\left[-\theta,0\right];\R^{N}\right)$.
The gain is given by $\gamma(r):=\max\{e^{d},1\}\max_{i}\tfrac{1}{s_i}\gamma_i(r).$
\end{theorem}

\textit{\textbf{Proof.}}
Let $0\neq x^t=((x^t_1)^T,\ldots,(x^t_n)^T)^T$ and $V$ be defined by $V(x^t):=\max_i\{\tfrac{1}{s_i}(V_i(x^t_i))\}$. For any $i\in\{1,\ldots,n\}$ consider open domains $M_i\in\R^N\backslash\{0\}$ defined by
\begin{align*}
M_i:=&\{\left((x^t_1)^T,\ldots,(x^t_n)^T\right)^T \in PC\left(\left[-\theta,0\right];\R^{N}\right):
%&\phantom{\{}
\tfrac{1}{s_i}V_i(x_i^t)>\max_{j\neq i}\tfrac{1}{s_j}V_j(x_j^t)\}.
\end{align*}

Take arbitrary $\hat{x}^t=((\hat{x}^t_1)^T,\ldots,(\hat{x}^t_n)^T)^T$ $\in PC\left(\left[-\theta,0\right];\R^{N}\right)$ for which there exist $i\in\{1,\ldots,n\}$, such that $\hat{x}^t\in M_i$. It follows that there is a neighborhood $U$ of $\hat{x}$ such that $V(x)=\frac{1}{s_i}V_i(x_i)$ and $D^+ V(x)$ exists.

%Now for any $\hat{x}^t=((\hat{x}^t_1)^T,\ldots,(\hat{x}^t_n)^T)^T$ $\in PC\left(\left[-\theta,0\right];\R^{N}\right)$ with\linebreak
%$(V_1(\hat{x}^t_1),\ldots,V_n(\hat{x}^t_n))^T\in M_i$ there is a neighborhood $U$ of $\hat{x}^t$ such that $V(x^t)=\frac{1}{s_i}V_i(x^t_i)$ and $D^+ V(x^t,u)$ exists.
By similar calculations as in the proof of Theorem~\ref{thm:main1}, $V$ is the exponential ISS-Lyapunov-Krasovskii functional for the overall system of the form \eqref{eq:SingleIMPTDSstatejump}. We can apply Theorem~\ref{th:SingleIMPTDSISSLyaKra} and the whole system is uniformly ISS over $\Simp [\mu,\lambda]$.\hfill $\Box$

%%%%%%%%%%%%%%%%%%%%%%%%%%%%%%%%%%%%%%%%%%%%%%%%%%%%%%%%%%%%%%%%%%%%%%%%%%%%%%%%%%%%%%%%%%%%%%%%%%%%%%%%%%%%%%%%%%%%%%%%%%%%%%%%%%%%%%%%%%%%%%%%%%%%%%%%%%%%%%%%%%%
\section{Example: networked control systems with time-delays}

We consider a class of networked control systems given by an interconnection of linear systems with time-delays \cite{HLT08}, \cite{WYB99}, \cite{NeT04} and \cite{HNX07}.  The $i$th subsystem is described as follows
\begin{equation}\label{ex:ncs_xi}
\begin{array}{lll}
\dot{x}_i&=&-a_{i}x_i+\sum\limits_{j,j\neq i}a_{ij}x_j(t-\tau_{ij})+b_i\nu_i,\ a_i>0,\\
y_i&=&x_i+\mu_i,\quad\quad i=1,\ldots,n.
\end{array}
\end{equation}
Here $\tau_{ij} \in [0,\theta]$ is a time-delay of the input from other subsystems with maximum involved delay $\theta>0$, $\nu_i$ is an input disturbance, $\mu_i$ a measurement/quantization noise. The sequence $\{t_1,t_2,\ldots\}$ is a sequence of time instances at which measurements of $x_i$ are sent. It is allowed to send only one measurement per each time instant. Between the sending of new measurements the estimate $\hat{x}_i$ of $x_i$ is given by
\begin{equation}\label{ex:ncs_hat{x}i}
\begin{array}{lll}
\dot{\hat{x}}_i(t)&=&-a_{i}\hat{x}_i(t)+\sum\limits_{j,j\neq i}a_{ij}\hat{x}_j(t-\tau_{ij}),\ t\not\in\{t_1,t_2,\ldots\}.\\
\end{array}
\end{equation}
At time $t_k$ the node $i_k$ gets access to the measurement $y_{i_k}$ of $x_{i_k}$ and all other nodes stay unchanged:
$$
\hat{x}_i(t_k)=\left\{
\begin{array}{ll}
y_{i_k}^{-}(t_k), & i=i_k,\\
\hat{x}_i^{-}(t_k), & i\neq i_k.
\end{array}\right.
$$
An estimation error is defined by $e_i:=\hat{x}_i-x_i$. The dynamics of $e_i$ can be then given by the following impulsive system:
\begin{equation}\label{ex:ncs_ei_de}
\begin{array}{lll}
\dot{e}_i(t)&=&-a_{i}e_i(t)+\sum\limits_{j,j\neq i}a_{ij}e_j(t-\tau_{ij})-b_i\nu_i,\ t\neq t_k,k\in\N,\\
\end{array}
\end{equation}
\begin{equation}\label{ex:ncs_ei_jump}
e_i(t_k)=\left\{
\begin{array}{ll}
\mu_{i_k}^{-}(t_k), & i=i_k,\\
e_i^{-}(t_k), & i\neq i_k.
\end{array}\right.
\end{equation}
The decision to which node a measurement will be sent is performed using some protocol, for examples see \cite{NeT04}.

%The TOD-like protocol sends a measurement to node $i_k$ if it has the largest absolute value error $\hat{x}_i^{-}-y_{i}^{-}=e_i^{-}-\mu_i^{-}$, i.e.
%\[|e_{i_k}^{-}-\mu_{i_k}^{-}|\geq |e_i^{-}-\mu_i^{-}|, i=1,\ldots,n.\]
Let us show that the error of the whole interconnected system \eqref{ex:ncs_ei_de}, \eqref{ex:ncs_ei_jump} is uniformly ISS using Lyapunov-Razumikhin approach. Firstly, we will find an ISS-Lyapunov-Razumikhin function candidate for each subsystem. 

Consider the function $V_i(e_i):=|e_i|$. If $t=t_k$, then $V_i(g_{i}(e_{i}))\leq\max\{|e_{i}|,|\mu_{i}|\}=\max\{e^{-d_{i}}V_i(e_{i}),|\mu_{i}|\}$ with $d_{i}=0$. Consider now the case $t\neq t_k$. If $|e_i|\geq \max\{\max\limits_{j,j\neq i}n\frac{|a_{ij}|}{a_i-\epsilon_i}  \max_{t-\theta \leq s \leq t}V_j(e_j(s)),n\frac{|b_i\nu_i|}{a_i-\epsilon_i}\}$, $\epsilon_i\in[0,a_i)$, then 
\begin{align*} 
	D^+V_i(e_i)=&(-a_{i}e_i+\sum\limits_{j,j\neq i}a_{ij}e_j(t-\tau_{ij})-b_i\nu_i)\cdot\sign e_i\\
\leq&-a_{i}|e_i|+\sum\limits_{j,j\neq i}|a_{ij}||e_j(t-\tau_{ij})|+|b_i\nu_i|\\
\leq&-a_{i}|e_i|+(a_i-\epsilon_i)|e_i|\\
=&-\epsilon_i|e_i|=-\epsilon_i V_i(e_i)=:-c_iV_i(e_i)
\end{align*}
%(a_i-\epsilon_i)e_i^2=-a_{i}e_i^2+a_ie_i^2-\epsilon_ie_i^2=-\epsilon_iV(e_i)=:-c_iV(e_i)$.
Thus, the function $V_i(e_i)=|e_i|$ is an exponential ISS-Lyapunov-Razumikhin function for the $i$th subsystem with $c_i=\epsilon_i$, $d_i=0$, $\gamma_{ij}(r)=n\frac{|a_{ij}|}{a_i-\epsilon_i}r$ and $\gamma_i(|(\mu_i,\nu_i)|)=\max\{1,n\frac{|b_i|}{a_i-\epsilon_i}\}|(\mu_i,\nu_i)|$.

To prove ISS of the whole error system we need to check the dwell-time condition \eqref{eq:cond_int} and the small-gain condition \eqref{eq:sgcequivalentLR}, see Theorem~\ref{th:main2}. Let us check condition \eqref{eq:cond_int}. We have $d=\min\limits_{i}d_i=0$, $c=\min\limits_{i}c_i=\min\limits_{i}\epsilon_i>0$. Taking $0<\lambda\leq c$ and any $\mu>0$ the dwell-time condition is satisfied for any $t\geq s\geq 0$ and time sequence $\{t_k\}$:
\[-dN(t,s)-(c-\lambda)(t-s)=-(c-\lambda)(t-s)\leq 0< \mu.\]
The fulfillment of the small gain condition \eqref{eq:sgcequivalentLR} can be checked by slightly modifying the cycle condition \cite{Rue07}: 
for all $(k_1,...,k_p) \in \{1,...,n\}^p$, where $k_1=k_p$, it holds
\begin{align}\notag
\gamma_{k_1k_2} \circ \gamma_{k_2k_3} \circ ... \circ \gamma_{k_{p-1}k_p} < e^{-\mu}\Id,
\end{align} 
where in this example we can choose $\mu$ arbitrarily small. If the cycle condition is fulfilled, then the the small gain condition is satisfied. Let us check this condition for the following parameters: 
$n=3$, %$\mu_1=0.01$, $\mu_2=0.03$, $\mu_3=0.02$, 
$b_i=1$, $\nu_i=2$, $\epsilon_i=0.1$, $i=1,2,3$; $\tau_{ij}=\theta=0.03$, $i,j=1,2,3$, $i \neq j$; $e(s)=(0.9; 0.3; 0.6)^T$, $s \in [-\theta,0]$;  $a_1=1$, $a_2=2$, $a_3=0.5$, 
\begin{align*}
	A:=(a_{ij})_{3\times 3}=\left(
\begin{array}{ccc}
0 & 0.25 & 0.25\\
0.7 & 0 & 0.65\\
0.15& 0.1&0
\end{array}\right).
\end{align*}
The system uses TOD-like protocol \cite{NeT04}. The protocol sends measurements at $t_k=0.1k,k\in\N$.

\begin{figure}[ht]
	\centering
		\includegraphics[width=0.9 \textwidth]{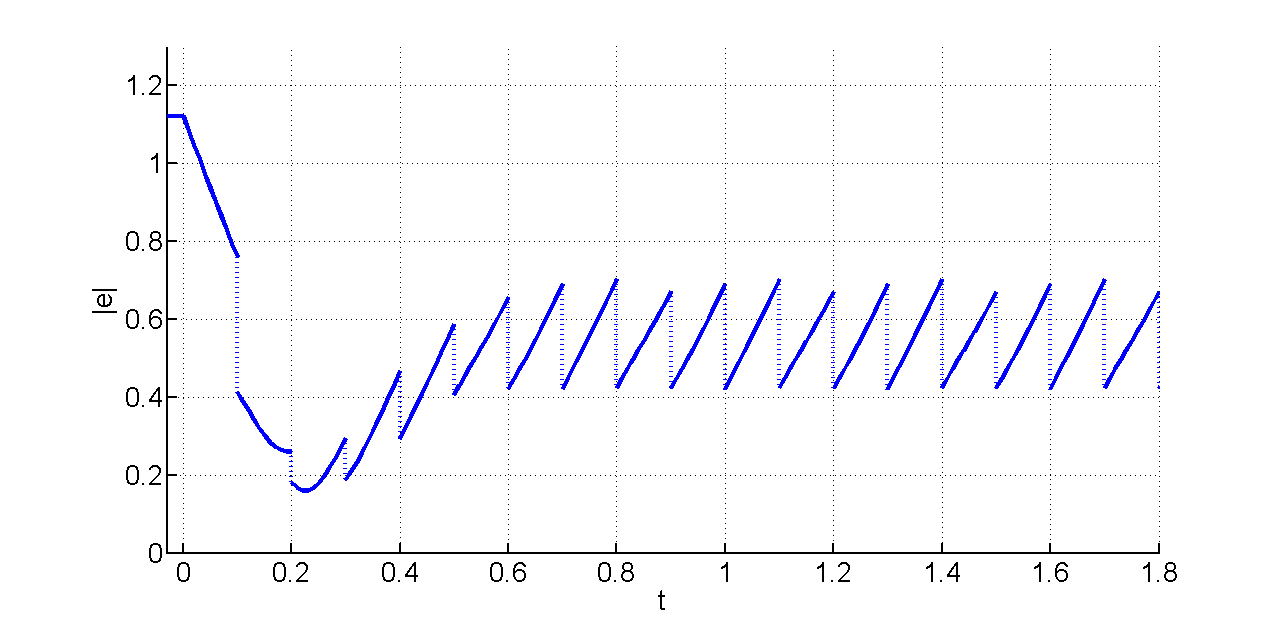}
		\caption{The trajectory of the Euclidean norm of the error of the networked control system.}
	\label{fig:NCS_example}
\end{figure}

The gain matrix $\Gamma$ is then given by
\begin{align*}
	\Gamma:=(\gamma_{ij})_{3\times 3}=\left(
\begin{array}{ccc}
0 & 0.8333 & 0.8333\\
1.1053   & 0 & 1.0263\\
1.1250 & 0.7500  &0
\end{array}\right).
\end{align*}
It is easy to check that all the cycles are less than the identity function multiplied by $e^{\mu}$, because $\mu$ can be chosen arbitrarily small. Thus, the cycle condition is satisfied and by application of Theorem~\ref{th:main2} the error system \eqref{ex:ncs_ei_de}, \eqref{ex:ncs_ei_jump} is uniformly ISS. The trajectory of the Euclidean norm of the error is given in Figure~\ref{fig:NCS_example}.

%%%%%%%%%%%%%%%%%%%%%%%%%%%%%%%%%%%%%%%%%%%%%%%%%%%%%%%%%%%%%%%%%%%%%%%%%%%%%%%%%%%%%%%%%%%%%%%%%%%%%%%%%%%%%%%%%%%%%%%%%%%%%%%%%%%%%%%%%%%%%%%%%%%%%%%%%%%%%%%%%%%
\section{Conclusions}
In this paper, we have established several theorems: At first we have introduced the Lyapunov-Krasovskii methodology and the Lyapunov-Razumikhin approach for establishing ISS of single impulsive systems with time-delays. Then, we have considered networks of impulsive subsystems without time-delays. As one of the results, we have proved an ISS-Lyapunov small-gain theorem, which guarantees that the whole network has the ISS property under a small-gain condition with linear gains and a dwell-time condition. To prove this, we have constructed the ISS-Lyapunov function and the gain of the whole system.
Under consideration of time-delays in such networks, we have proved two theorems to show that a network has the ISS property provided that a small-gain condition with linear gains and a dwell-time condition is satisfied. On the one hand, we have used ISS-Lyapunov-Razumikhin functions and on the other hand we have used ISS-Lyapunov-Krasovskii functionals. An application was illustrated for networked control systems with time delays. 

An open question to be investigated in the future is the usage of general Lyapunov functions instead of exponential Lyapunov functions for single systems, where a different dwell-time condition has to be formulated. 
Then, for the interconnection one can also use general Lyapunov functions and general gains instead of only linear gains. 
Furthermore, for interconnected systems the case could be investigated the 
case, when the impulse sequences of subsystems are different.

%The relation between the small-gain condition and the dwell-time condition should be investigated.

\bibliographystyle{elsarticle-num}
%\bibliography{literatur}

\end{document}